\documentclass[11pt]{article}
\usepackage{geometry}
\usepackage{graphicx} 
\usepackage{amsmath}
\usepackage{amsthm}
\usepackage[toc,page]{appendix}
\usepackage{tabularray}
\usepackage{amssymb}
\usepackage{caption}
\usepackage{mathrsfs}
\usepackage{eufrak, color}
\usepackage{indentfirst}
\usepackage{array}
\usepackage{tikz}
\usepackage{enumerate}
\usepackage[shortlabels]{enumitem}
\usetikzlibrary{graphs}
\usetikzlibrary{graphs.standard}
\usepackage{hyperref}
\newcommand{\setEnvironmentQed}[2]{
  \AtBeginEnvironment{#1}{%
    \pushQED{\qed}\renewcommand{\qedsymbol}{#2}%
  }
  \AtEndEnvironment{#1}{\popQED}
}
\newtheorem*{example}{Example}
\newtheorem{lemma}{Lemma}
\newtheorem{theorem}{Theorem}
\setEnvironmentQed{example}{\ensuremath{\blacksquare}} 
\usepackage{doi}
\usepackage{url}
\usepackage{comment}
\usepackage{enumitem}
\providecommand{\keywords}[1]
{
  \small	
  Keywords:  #1
}

\title{\Large{Exact upper bounds for the minimum sizes\\ of strong and weak separating path systems of cliques}}
\author{George Kontogeorgiou\thanks{Center for Mathematical Modeling (CNRS IRL2807), University of Chile. Supported by ANID Basal Grant CMM FB210005 and ANID-FONDECYT Postdoctorado Grant No. 3250479. {\it Email:} gkontogeorgiou@dim.uchile.cl
}
~and Maya Stein\thanks{Department of Mathematical Engineering and Center for Mathematical Modeling (CNRS IRL2807), University of Chile. Supported by FONDECYT Regular Grant 1221905,  by ANID Basal Grant CMM FB210005,  and by MSCA-RISE-2020-101007705 project {\it RandNET}.
{\it Email:} mstein@dim.uchile.cl
}}

\begin{document}

\maketitle

\begin{abstract}
We prove an upper bound of $n+9$ for the strong  separation number of the complete graph $K_n$, and an upper bound of $n+1$ for its weak separation number. This improves on the previous best known bound of $(1+o(1))n$ for both cases.
\end{abstract}

\vspace{5mm}

\section{Introduction}

Let $G$ be a graph. A \emph{path system} of $G$ is a set of paths in $G$. Given edges $e,f\in E(G)$ and a path system~$\mathcal{P}$ of $G$, we say that $e$ is \emph{separated from} $f$ by $\mathcal{P}$ if there exists a path $P\in\mathcal{P}$ such that $e\in P$ and $f\notin P$. A path system $\mathcal{P}$ \emph{weakly separates} $G$ if for every pair of edges of $G$ it separates one of them from the other, and \emph{strongly separates} $G$ if it separates all the edges of $G$ from each other. Let $ssp(G)$ (resp.~$wsp(G)$) be the minimum size of a path system that separates $G$ strongly (resp.~weakly), called the \emph{strong} (resp. \emph{weak}) \emph{separation number} of $G$. Clearly, $wsp(G)\leq ssp(G)$ for every graph $G$. We are interested in the problem of determining these two quantities when $G$ is a complete graph on $n$ vertices.

The study of general separating set systems was initiated by R\'enyi in the 1960s \cite{renyi}. Various versions concerning the separation of edges by subgraphs were introduced around the 2000s by computer scientists who were seeking efficient ways to detect faulty links in 
networks~\cite{compsci1, compsci2, compsci3, compsci4}.
The problem of determining tight upper bounds for the size of separating path systems of graphs  on $n$ vertices 
has become very popular in the combinatorics community during the last decade. 

Among the first to consider this problem were Balogh, Csaba, Martin and Pluh\'ar in the strong setting~\cite{balogh2016path}, and Falgas-Ravry, Kittipassorn, Kor\'andi, Letzter and Narayanan~\cite{nlogn} in the weak setting. Both groups proved a common upper bound of $O(n \log n)$ and conjectured one of $O(n)$. Letzter \cite{nlog*n} came very close to the conjectured bound by proving that $ssp(G)=O(n \log^*n)$. The conjecture was eventually settled by Bonamy, Botler, Dross, Naia and Skokan~\cite{19n}, who proved that $ssp(G)\leq 19n$. They 
further conjectured that $(1+o(1))n$ is optimal.

Restricting the question to the class of complete graphs is an interesting and appropriate direction to follow, as noted e.g.~by Letzter in the conclusion of \cite{nlog*n}. This is not only because of the reduced complexity of the problem but also because of the fact that complete graphs have the highest known lower bound for their weak separation number among all graph classes. It is therefore natural to suspect that they are extremal for this problem. Specifically, the lower bound $wsp(K_n)\geq n-1$ is tight \cite{nlogn}. As for the upper bound, owing to work by Fernandes, Oliveira Mota and Sanhueza-Matamala \cite{Nico}, it is known that \[ssp(K_n)=(1+o(1))n.\]

Wickes \cite{wickes} showed that every complete graph of odd prime order has a generating path, and that the set of all the images under rotation of that generating path is a weakly separating path system. So, $wsp(K_n)\leq n$ when $n$ is prime. She also proved that $wsp(K_n)\leq n$ when $n-1$ is prime or when $n\leq 20$. 

In Section \ref{sec3}, we give a similar upper bound for $wsp(K_n)$ for all $n\in\mathbb N$. 

\begin{theorem} \label{theorem1}
    For every natural number $n$, $wsp(K_n)\leq n+1$.
\end{theorem}

In combination with the lower bound of $n-1$, our result nearly resolves the problem of the size of minimum weakly separating path systems for complete graphs.
Moreover, our proof is both short and entirely constructive, providing concrete and relatively easy to describe examples of weakly separating path systems for all complete graphs. 

As we explain in Section \ref{sec4}, the path system constructed in the proof of Theorem \ref{theorem1} leaves few edges that are not separated from every other edge, and we will show that these may be assumed to be spaced in a relatively convenient way. Consequently, we show that we can expand the path system of Theorem \ref{theorem1} by eight more paths to make it strongly separating.

\begin{theorem} \label{theorem2}
    For every natural number $n$, $ssp(K_n)\leq n+9$.
\end{theorem}

\section{Definitions and Preliminaries}

 A \emph{labelled clique} is a complete graph $K_n$ together with a labelling of $V(K_n)$ with numbers from an interval of length $n$ of the natural numbers. The labelling naturally induces on $V(K_n)$ a cyclic ordering, which we interpret as counterclockwise, and a rotation action. Given a labelled clique $K_n$, we partition $E(K_n)$ in \emph{types} as follows: 
   the type of an edge $\{x_1, x_2\}$ is defined to be $\min\{(x_1-x_2) 
 (\text{mod $n$}), (x_2-x_1) 
 (\text{mod $n$})\}$. The type  of every edge is in $[\lceil\frac{n-1}{2}\rceil]$. 
  Edges of the same type are called \emph{homotypical}. 
  
  A \emph{counterclockwise rotation} of order $i$ in a labelled clique $K_n$ is an automorphism that maps the vertex $\ell$ to the vertex $(\ell +i)\text{ (mod $n$)}$ for each $\ell\in [n]$.
  The \emph{counterclockwise distance} of two homotypical edges is the minimum order of a counterclockwise rotation that maps one to the other. Equivalently, for homotypical edges $\{x_1, x_2\}$, $\{y_1, y_2\}$ with, say, $(x_1-x_2) 
 \text{ (mod $n$)}=(y_1-y_2) 
 \text{ (mod $n$)}$, their counterclockwise distance is the type of $\{x_1,y_1\}$.

 We define a \emph{generating path} $P$ of $K_n$ as one that exhibits the following properties: 

\begin{itemize}[leftmargin=13mm]
   \item[\textbf{(GP1)}] 
    there is a type $c$ such that $P$ contains exactly one edge of type $c$;
    \item[\textbf{(GP2)}] $P$ contains exactly two edges of each other type;
    \item[\textbf{(GP3)}] any two pairs of homotypical edges in $P$ have distinct counterclockwise distances.
\end{itemize}

The \emph{set of images under rotation of a   path~$P$} is the family of all paths obtained by rotating~$P$ counter\-clockwise by $0,1,\ldots,n-1$. Note that rotation maps any edge of $P$ to an edge of the same type. 
Now let $P$ be a generating path, and let $c$ be as in  {\textbf{(GP1)}. Then, if $n$ is odd,  properties {\textbf{(GP1)} and {\textbf{(GP2)} imply that,  in the set of images under rotation of $P$, each edge of type $c$ is in exactly one path, and each other edge is in exactly two paths}}}\footnote{This fails if $n$ is even, as then each edge of type $\frac{n}{2}$ lies in two paths if $c=\frac{n}{2}$, and in four paths otherwise.} 

The following statements are proved in \cite{wickes}. 

\begin{itemize}[leftmargin=11mm]

\item [\textbf{(W1)}] Every complete graph of prime order at least $5$ contains a \emph{generating cycle}, that is,~a cycle containing exactly two edges of each type and satisfying \textbf{(GP3)}. In order to obtain a generating path from a generating cycle, we simply delete a single edge.\footnote{In fact, Wickes \cite{wickes} proved that every complete graph of odd prime order contains a generating path. However, \textbf{(W1)} follows trivially from her proof, see Lemma \ref{generating_cycle} in Appendix \ref{appendix2}.} 

\item [\textbf{(W2)}] In every complete graph of odd prime order, the set of images under rotation of a generating path is a weakly separating path system. 

\item [\textbf{(W3)}] Every complete graph of order between $3$ and $19$ contains a path (not necessarily a generating path) whose images under rotation form a weakly separating path system in which only the edges of at most one type $c$ fail to be separated from all of the other edges. These paths and the corresponding types are listed in Lemma \ref{small_paths} of Appendix \ref{appendix2}.

\end{itemize}


\begin{example} In Figure \ref{fig:1} we see a generating cycle $C_{example}$ for $K_{17}$. Indeed, inspection reveals that the depicted cycle has, for each type in $\{1,\dots, 8\}$, two edges of that type. To verify $\textbf{(GP3)}$, we direct the reader to the rightmost column of Table \ref{tab1}. Note that $C_{example}$ spans its ambient clique but for one vertex. This is the case in general for the generating paths and cycles obtained in \cite{wickes}. We will not use this property, but we state it in order to aid the reader's intuition.

\end{example}

\begin{figure}[!h]
        \centering
        \caption{The generating cycle $C_{example}$.}
        \label{fig:1}
            \begin {tikzpicture}[-latex ,auto ,node distance =2cm and 3cm ,semithick , state/.style ={draw, circle}]
    \node[state,scale=0.9] (1) at (8,8) {1};
    \node[state,scale=0.9] (2) at (6.5,7.5) {2};
    \node[state,scale=0.9] (3) at (5.2,6.8) {3};
    \node[state,scale=0.9] (4) at (4.5,5.5) {4};
    \node[state,scale=0.9] (5) at (4.2,4) {5};
    \node[state,scale=0.9] (6) at (4.5,2.5) {6};
    \node[state,scale=0.9] (7) at (5.2,1.2){7};
    \node[state,scale=0.9] (8) at (6.5,0.5) {8};
    \node[state,scale=0.9] (9) at (8,0) {9};
    \node[state,scale=0.8] (10) at (9.5,0) {10};
    \node[state,scale=0.8] (11) at (11,0.5) {11};
    \node[state,scale=0.8] (12) at (12.3,1.2){12};
    \node[state,scale=0.8] (13) at (13,2.5) {13};
    \node[state,scale=0.8] (14) at (13.2,4.2) {14};
    \node[state,scale=0.8] (15) at (12.8,6) {15};
    \node[state,scale=0.8] (16) at (11.5,7.5) {16};
    \node[state,scale=0.8] (17) at (10,8) {17};

     \path[-] (17) edge[bend left=30] node {} (1);
     \path[-] (1) edge[bend left=30] node {} (4);
     \path[-]  (4) edge  node {} (13);
     \path[-] (13) edge  node {} (6);
     \path[-] (6) edge[bend right=30]  node {} (2);
     \path[-]  (2) edge[bend left=30]  node {} (7);
     \path[-]  (7) edge[bend right=30]  node {} (5);
     \path[-] (5) edge[bend right=30]  node {} (16);
     \path[-] (16) edge[bend right=30]  node {} (15);
     \path[-] (15) edge[bend right=30]  node {} (12); 
     \path[-] (12) edge  node {} (3); 
     \path[-] (3) edge  node {} (10);
     \path[-] (10) edge[bend left=30]  node {} (14);
     \path[-] (14) edge[bend right=30]  node {} (9);
     \path[-] (9) edge[bend left=30]  node {} (11);
     \path[-] (11) edge[bend left=30]  node {} (17);
    \end{tikzpicture}
\end{figure}

\begin{table}[ht]
 \begin{center}
 \caption{Its edges by type, and the counterclockwise distances of its pairs of homotypical edges.}
 \label{tab1}
\begin{tabular}{ c | c | c }
type & edges & counterclockwise distance\\
\hline
1 & $\{15,16\}$, $\{17,1\}$ & 2\\
2 & $\{5,7\}$, $\{9,11\}$ & 4\\
3 & $\{12,15\}$, $\{1,4\}$ & 6\\
4 & $\{2,6\}, \{10,14\}$ & 8\\
5 & $\{2,7\}$, $\{9,14\}$ & 7\\
6 & $\{11,17\}$, $\{16,5\}$ & 5\\
7 & $\{3,10\}$, $\{6,13\}$ & 3\\
8 & $\{12,3\}$, $\{13,4\}$ & 1
\end{tabular}
\label{typetable}
\end{center}
\end{table}

Finally, we work under the convention that every path is \emph{oriented}, in the sense that it has has a first and a last vertex. If the last vertex $u$ of a path $P$ is distinct from the first vertex $v$ of a path $P'$, then
we denote by $e(P,P')$ the edge connecting $u$ and $v$.

\subsection{Overview of the Construction}

Before fleshing out the details, we provide a bird's eye view of our work.

\subsubsection{Weakly separating systems: Proof sketch for Theorem 1 (Section \ref{sec3})}

In the beginning, we decompose $n$ into a rapidly decreasing sequence of odd prime numbers $p_1,\dots,p_m$ and a small rest $b\in\{3,\dots,19\}$. Then, we decompose $K_n$ into vertex-disjoint cliques $K_{p_1},\dots,K_{p_m}$ and $K_b$. In each of the cliques $K_{p_k}$, we find a generating cycle $C_k^1$.  These are known to exist by~\textbf{(W1)}. For each $k\in [m]$, we denote by $\mathcal{C}_k$ the family of images under rotation of $C_k^1$ in $K_{p_k}$. In particular, $|\mathcal{C}_k|=p_k$ for each $k\in [m]$. We also find in $K_b$ a path $\mathbf{B}_{m+1}^1$ such that the set $\mathscr{B}$ of its images under rotation in $K_b$ is a weakly separating path system of $K_b$. We know that such a path exists by \textbf{(W3)}. Note that $|\mathbf{B}_{m+1}|=b$.

For each $k\in [m]$, we can choose a type $c_k$, delete an edge of this type from $C_k^1$, and then delete the images under rotation of said edge from the corresponding cycles in $\mathcal{C}_k$ to obtain a family of paths $\mathcal{P}_k$ that weakly separates the edges of $K_{p_k}$ by \textbf{(W2)}. Then the family $\left(\bigcup_{k=1}^m\mathcal{P}_k\right)\cup \mathscr{B}$ has size $n$ and weakly separates the set of edges $\left(\bigcup_{k=1}^mE(K_{p_k})\right)\cup E(K_b)$. However, this fact is not 
sufficient for
us. The problem is that the edges that go \emph{between} the cliques $K_{p_1},\dots , K_{p_m}, K_b$ would remain unseparated (and even uncovered). To remedy this, we must firstly modify our cycle systems.

Let us visualise the set of clique components $K_{p_1},\dots, K_{p_m},K_b$ as a tower, with each successive floor being much smaller than the floors below it. We will call these floors \emph{levels}. Each member of $\left(\bigcup_{k=1}^m\mathcal{C}_k\right)\cup\mathscr{B}$ lives in one of the levels. Let us focus on the level of $K_{p_k}$. A good way to fit into our cycles some of the edges that go in-between distinct floors is to take $C_k^1$, choose an edge $e_1$ of some type $t_1$ inside it, and substitute $e_1$ by a path of length $2$ 
going through a vertex $v_1$ that lives in a higher level. Then we do the same for some other edge $e_2$ of a different type $t_2$ in $C_k^1$ and for a different vertex $v_2$ at a higher level, and so on. Because the sequence $p_1,\dots, p_m,b$ decreases so rapidly that the size of $K_{p_k}$ is much larger than the total size of all the levels above it, there are enough edges of distinct types in $C_k^1$ to do this modification for every higher vertex. 
Finally, we rotate the modified $C_k^1$ except for its new vertices $v_i$ which stay fixed (it may help to think of the edges of the modified $C^1_k$ incident to $v_i$ as rubber bands). In this way we obtain a new cycle system, whose members jointly contain all the edges that either live in $K_{p_k}$, or have one end in $K_{p_k}$ and one end at a higher level. There are two problems with this approach 
which we will discuss separately.

\paragraph{I. The new edges are not separated from a specific edge of the cycle to which they were added.}
Let $a_1$ and $b_1$ be the edges of the path of length $2$ that substituted $e_1$ in the above process, that is, the one that 
goes through
$v_1$. As the type of $e_1$ is $t_1$, we observe that counterclockwise $t_1$-rotation maps one of $a_1$, $b_1$ to the other, let us say the former to the latter. 
The same is true for a pair of homotypical edges in $K_{p_k}$
that have counterclockwise distance $t_1$ --- exactly one such pair exists by the definition of a generating cycle. 
Let $e$ and $e'$ be two such edges, with $e$ being mapped to $e'$ by counterclockwise $t_1$-rotation. Then
$e'$ and $b_1$ 
both lie in 
the modified cycle $C_k^1$ and its image under counterclockwise $t_1$-rotation. They lie in no other cycles.
Therefore, the edges $b_1$ and $e'$ are not weakly separated from each other by $\mathcal{C}_k$. The resolution to this sort of issue is an important part of our proof, and
 we will discuss it now extensively.    

We begin by asking: what if such an inconvenient homotypical pair $e$, $e'$ in $K_{p_k}$ did not exist? We actually get a say over this. When we choose the next edge to be substituted, we may very well choose $e_2:=e'$ or $e_2:=e$. Now, $b_1$ will be separated from $e'$, as the latter does not belong to the modified $C_k^1$ (resp. to its image under counterclockwise $t_1$-rotation). Of course, we then  have the same problem with an edge incident to $v_2$ and an edge in a homotypical pair of counterclockwise distance $t_2$, so we pick one of the edges of said pair to be $e_3$, and so on. Each time, we choose the next edge to be substituted from the homotypical pair of edges whose counterclockwise distance is equal to the type of the last edge that was substituted. If at some point it is impossible to choose the next edge to be substituted like this (because there does not remain such a pair, that is, the type of this pair has already been used once), then we arbitrarily choose an edge of a type that is yet unused. We do all this in every cycle of $\mathcal{C}_k$, in a manner consistent under rotation.   

But then what happens at the last step? Once again, there is an edge, between a vertex of $K_{p_k}$ and a higher vertex, that is not separated from some edge in $K_{p_k}$. But now there are no more higher vertices to span, so if we delete this edge of $K_{p_k}$ from $C_k^1$, we cannot ``bridge the gap" with a path of length $2$. So, the modified $C_k^1$ will not be a cycle at all, but a path, and its images under rotation will form a path system. But this is exactly what we want! In fact, the reason that we started by taking a generating cycle in each clique (instead of a generating path) was so that, when we would eventually be forced to simply delete an edge, we would be glad about it. In the proof, we denote the type of that last deleted edge by $c_k$, and we actually compute it and remove it first rather than at the end.

In this manner, for each $k\in [m]$, we obtain a path system $\mathscr{P}_k$ of size $p_k$ that separates a large proportion of the edges that have their lower end in $K_{p_k}$. We do not change $\mathscr{B}$.

\paragraph{II. The substituted and the deleted edges from the cycles are no longer separated.} Many of the edges of $K_{p_k}$ were deleted from the cycles of $\mathcal{C}_k$ to form the paths of $\mathscr{P}_k$. As a result, they are not in enough of these paths and they are not separated from many other edges. However, we can gather most of them in a path system $\mathscr{R}_k$. Each path in this system contains every edge (except one) of a type from which we deleted edges, as well as one carefully chosen extremal edge that connects it to the next higher level. In particular, the size of $\mathscr{R}_k$ is one more than the number of vertices in levels higher than $K_{p_k}$ (because, apart from the type $c_k$, every other type from which we deleted an edge witnesses the existence of a higher vertex through which we reconnected, with a path of length $2$, the ends of said deleted edge).

We now have $n$ paths in $(\bigcup_{k=1}^m\mathscr{P}_k)\cup\mathscr{B}$, each starting at a different vertex of $K_n$ and living in its corresponding level and the levels above, and a lot more paths in $\bigcup_{i=1}^m\mathscr{R}_k$, each contained in its level but for one extremal edge. These extremal edges have been chosen in a way that allows certain concatenations: essentially, to each path in $(\bigcup_{k=1}^m\mathscr{P}_k)\cup\mathscr{B}$ we attach a ``tail" that consists of many paths in $\bigcup_{i=1}^m\mathscr{R}_k$, one from each lower level. We also concatenate all of the paths in $\bigcup_{i=1}^m\mathscr{R}_k$ that consist of edges of the special type $c_k$ for each $k\in [m]$; there is one such at each level up to $m$. After doing these concatenations, we are left with a family of paths $\mathscr{Q}$, whose members are the paths in $(\bigcup_{k=1}^m\mathscr{P}_k)\cup\mathscr{B}$ modified with the addition of $\mathscr{R}$-tails, as well as one more path that collects edges of the special type $c_k$ from each level $K_{p_k}$. Therefore, $|\mathscr{Q}|=|(\bigcup_{k=1}^m\mathscr{P}_k)\cup\mathscr{B}|+1=n+1$.

Lastly, we 
verify by case analysis that $\mathscr{Q}$ is indeed a weakly separating path system for $K_n$.

\subsubsection{Strongly separating systems: Proof sketch for Theorem 2 (Section \ref{sec4})}

We proceed to adapt these ideas to the strong setting. To begin, we need a decomposition of $K_n$ in which the sizes of our tower levels decrease even faster than before. This is always possible but for a small number of exceptions that we list and consider separately when necessary. Then we apply anew our previous construction to obtain a weakly separating path system of $K_n$ based on this new decomposition.   

We point out that our weakly separating path system is already almost a strongly separating path system; there is only a small number of edges that are not separated from all of the other edges. Our solution for this is as follows. In Section \ref{sec3}, when first designing our weakly separating path system, we choose each family of paths $\mathscr{R}_k$ almost arbitrarily, with only one easy-to-fulfill condition to guide our choice. 
Then, in Section \ref{sec4},
we make a very careful choice obeying a number of conditions, which helps us ensure that the offending edges all lie inside at most six linear forests. It is then a simple task to separate them from all of the other edges by using only eight more paths, resulting in the final bound of $n+9$ for the strong separation number of a complete graph on $n$ vertices.   \\

\section{Proof that $wsp(K_n)\leq n+1$}\label{sec3}

Let $3\leq n=\sum_{k=1}^mp_k+b$ with $p_k$ odd primes. For brevity, set $s_0:=0$ and $s_k:=p_1+\dots+p_k$ for $k\in [m]$. It is a corollary of \cite{nagura} (see Lemma \ref{lem1} in Appendix \ref{appendix1}) that we can choose each $p_k$ so that $n-s_k\leq \frac{p_k-3}{2}$ and $b\in \{3,\dots,19\}$. In particular, $p_k>p_{k+1}$ for $1\le k<m$, and $p_m>b\geq 3$, so $p_m\geq 5$.

We define pairwise vertex-disjoint induced cliques of $K_n$ indexed by, and with sizes corresponding to, the terms of the above sum. We then label the vertices of each clique $K_{p_k}$ using the interval $[s_{k-1}+1, s_k]$. We call these labelled cliques \emph{levels}, and we say that an edge \emph{touches} a level if its smaller end is in that level.
As each level is a clique, if an edge is in a level, then we can
think of its type with respect to that clique, that is, with respect
to that level. In what follows, we only refer to types of edges that
are contained in a level, and the type then refers to that level.

\begin{example}[cont]

 Consider $K_{21}$. It can be decomposed into two vertex-disjoint cliques, $K_{17}$ and $K_4$, for which our size requirements are satisfied. That is, $17$ is an odd prime, $3\leq 4\leq 19$, and $21-17\leq\frac{17-3}{2}$. We label the vertices of $K_{17}$ with $\{1,\dots, 17\}$, and those of $K_4$ with $\{18,\dots,21\}$.
    
\end{example}

Let \[\mathscr{B}:=\{\mathbf{B_{m+1}^{s_m+1}},\dots,\mathbf{B_{m+1}^n}\}\] be a weakly separating path system for $K_b$ given by the images under rotation of the corresponding path from \textbf{(W3)}. Each exponent denotes the first vertex of the corresponding path.

\begin{example}[cont]
    A weakly separating path system of size $4$ for our vertex-labelled $K_4$ is the one generated by rotating the path $[18,19,21]$. Explicitly, we have $\mathscr{B}:=\{[18,19,21],[19,20,18],[20,21,19],[21,18,20]\}$. 
\end{example}

For each $k\in [m]$, let $C_k^{s_{k-1}+1}$ be a generating cycle for the clique $K_{p_k}$, which exists by \textbf{(W1)}. Let \[\mathcal{C}_k:=\{C_k^{s_{k-1}+1},\dots,C_k^{s_k}\}\] be the set of its images under rotation, with $C_k^{s_{k-1}+i}$ being the image of $C_k^{s_{k-1}+1}$ under counterclockwise $(i-1)$-rotation. 
We construct a weakly separating path system $\mathscr{Q}$ for $K_n$ as follows.

To begin, we note that the generating cycle $C_k^{s_{k-1}+1}$ defines a permutation $\sigma_k\in S_{\frac{p_k-1}{2}}$ that maps each type in $[\frac{p_k-1}{2}]$ to the counterclockwise distance in $K_{p_k}$ of the homotypical edges of $C_k^{s_{k-1}+1}$ of that type. Indeed, the map $\sigma_k$ is a permutation of $[\frac{p_k-1}{2}]$ because it is an injection (by \textbf{(GP3)}), hence a bijection (because it is a function between two equinumerous finite sets). Based on this observation, for each $p_k$ we define an injection $\phi_k:\{s_k+1,\dots,n\}\rightarrow [\frac{p_k-1}{2}]$ in the following way: 
$\phi_k(s_k+1)$ is arbitrary, and for each $j\in\{2,\dots,n\}$, we choose $\phi_k(s_k+j)$ as the counterclockwise distance of the two edges of type $\phi_k(s_k+j-1)$ in $C_k^{s_{k-1}+1}$ (that is, $\phi_k(s_k+j):=\sigma_k\circ\phi_k(s_k+j-1)$) if this allows $\phi_k$ to retain its status as an injection; otherwise we choose $\phi_k(s_k+j)$ arbitrarily among the unused distances.
An equivalent way to think about it is that we determine the values of $\phi_k(s_k+1),\dots,\phi_k(n)$ by choosing a type in $[\frac{p_k-1}{2}]$ to be the value of $\phi_k(s_k+1)$, and then following its cycle in $\sigma_k$. When we complete this cycle, we choose a new type outside of it and follow the cycle of that, and so on. Note also that $\phi_k$ is not surjective, as $n-s_k\leq\frac{p_k-3}{2}$.


 We then select a special type $c_k\in [\frac{p_k-1}{2}]\setminus im(\phi_k)$. If $im(\phi_k)$ is a union of some cycles of $\sigma_k$, we let $c_k$ be arbitrary. Otherwise, we choose $c_k$ to be the counterclockwise distance of the edges of type $\phi_k(n)$. Then we pick an edge of type $c_k$ in $C_k^{s_{k-1}+1}$. For each $i\in [p_k]$, we remove from the cycle $C_k^{s_{k-1}+i}$ 
the image of said edge under counterclockwise $(i-1)$-rotation (which also has type $c_k$).
 This yields a generating path $P_k^{s_{k-1}+1}$ and the set of its images under rotation (which is a weakly separating path system of $K_{p_k}$ by \textbf{(W2)}), \[\mathcal{P}_k:=\{P_k^{s_{k-1}+1},\dots,P_k^{s_k}\}.\] Finally, we shift  
 the exponents of the paths of $\mathcal{P}_k$ so that the path $P_k^{s_{k-1}+j}$ starts from $s_{k-1}+j$. 
 
 Before we continue, as a useful reminder and an intuition check, let us recount some simple facts about the paths in $\mathcal{P}_k$: each path $P_k^{s_{k-1}+j}$ lies in level $k$ (inherited from $C_k^{s_{k-1}+j}$), starts at the vertex $s_{k-1}+j$\\ (we shifted the exponents to ensure this), ends at the vertex $s_{k-1}+\left(j-c_k \text{ (mod } p_k)\right)$
 (because together with an edge of type $c_k$ it forms the cycle $C_k^{s_{k-1}+j}$), and spans level $k$ but for one vertex (inherited from $C_k^{s_{k-1}+j}$).   

 \begin{example}[cont]

Consider the generating cycle $C_1^1:=C_{\text{example}}$ of $K_{17}$ from Figure \ref{fig:1}. Let us compute in our example of $K_{21}$, in the manner described above, the permutation $\sigma_1$ induced by $C_1^1$, as well as an appropriate injection $\phi_1$ and type $c_1$. To find $\sigma_1$, let us look at Table \ref{typetable}. We see that the edges of type $1$ have counterclockwise distance $2$, the edges of type $2$ have counterclockwise distance $4$, the edges of type $4$ have counterclockwise distance $8$, etc. In this way, we find that $\sigma_1=(1248)(3657)$. The value of $\phi_1(18)$ can be selected arbitrarily, let us say $\phi_1(18):=3$, but this choice determines also the other values of $\phi_1$: we have $\phi_1(19)=\sigma_1\circ\phi_1(18)=6$, $\phi_1(20)=\sigma_1\circ\phi_1(19)=5$ and $\phi_1(21)=\sigma_1\circ\phi_1(20)=7$. As $im(\phi_1)$ consists of the elements of a cycle of $\sigma_1$, we must choose the value of $c_1$ arbitrarily from $\{1,\dots,8\}\setminus im(\phi)=\{1,2,4,8\}$. We choose $c_1:=1$. For ease of inspection, we collect these values in Table \ref{tablevalues}. 

\begin{table}[ht]
\caption{The values of $\sigma_1$, $\phi_1$ and $c_1$ for our example on $K_{21}$.}
\label{tablevalues}
 \begin{center}
\begin{tabular}{ c | c | c }
$\sigma_1$ & $\phi_1$ & $c_1$ \\ 
 \hline
$(1248)(3657)$ & $\phi_1(18)=3$, $\phi_1(19)=6$ & 1 \\
  & $\phi_1(20)=5$, $\phi_1(21)=7$ & 
\end{tabular}
\end{center}
\end{table}

To obtain $P_1^1$ (shown in Figure \ref{fig:2}), we must choose and remove from $C_1^1$ an edge of type $c_1=1$. We choose the one incident to the vertex $1$, so that $P_1^1$ already starts from $1$ and no shifting of exponents is required. Each path $P_1^i\in\mathcal{P}_1$ is the image of $P_1^1$ under counterclockwise $(i-1)$-rotation. 

\end{example}

\begin{figure}[!h]
        \centering
        \caption{The generating path $P_1^1$ obtained from the generating cycle $C_1^1:=C_{example}$.}
        \label{fig:2}
        $\displaystyle
            \begin {tikzpicture}[-latex ,auto ,node distance =2cm and 3cm ,semithick , state/.style ={draw, circle}]
    \node[state,scale=0.9] (1) at (8,8) {1};
    \node[state,scale=0.9] (2) at (6.5,7.5) {2};
    \node[state,scale=0.9] (3) at (5.2,6.8) {3};
    \node[state,scale=0.9] (4) at (4.5,5.5) {4};
    \node[state,scale=0.9] (5) at (4.2,4) {5};
    \node[state,scale=0.9] (6) at (4.5,2.5) {6};
    \node[state,scale=0.9] (7) at (5.2,1.2){7};
    \node[state,scale=0.9] (8) at (6.5,0.5) {8};
    \node[state,scale=0.9] (9) at (8,0) {9};
    \node[state,scale=0.8] (10) at (9.5,0) {10};
    \node[state,scale=0.8] (11) at (11,0.5) {11};
    \node[state,scale=0.8] (12) at (12.3,1.2){12};
    \node[state,scale=0.8] (13) at (13,2.5) {13};
    \node[state,scale=0.8] (14) at (13.2,4.2) {14};
    \node[state,scale=0.8] (15) at (12.8,6) {15};
    \node[state,scale=0.8] (16) at (11.5,7.5) {16};
    \node[state,scale=0.8] (17) at (10,8) {17};

     \path[-] (1) edge[bend left=30] node {} (4);
     \path[-]  (4) edge  node {} (13);
     \path[-] (13) edge  node {} (6);
     \path[-] (6) edge[bend right=30]  node {} (2);
     \path[-]  (2) edge[bend left=30]  node {} (7);
     \path[-]  (7) edge[bend right=30]  node {} (5);
     \path[-] (5) edge[bend right=30]  node {} (16);
     \path[-] (16) edge[bend right=30]  node {} (15);
     \path[-] (15) edge[bend right=30]  node {} (12); 
     \path[-] (12) edge  node {} (3); 
     \path[-] (3) edge  node {} (10);
     \path[-] (10) edge[bend left=30]  node {} (14);
     \path[-] (14) edge[bend right=30]  node {} (9);
     \path[-] (9) edge[bend left=30]  node {} (11);
     \path[-] (11) edge[bend left=30]  node {} (17);
    \end{tikzpicture}
    $
    
\end{figure}
 
 Let us now modify the paths of each family $\mathcal{P}_k$ to also cover the edges that are not contained in any level. For each $j\in [n-s_k]$, we pick exactly one edge of type $\phi_k(s_k+j)$ in $P_k^{s_{k-1}+1}$, and for each $i\in [p_k]$, we remove from the path $P_k^{s_{k-1}+i}\in\mathcal{P}_k$ the image of said edge under counterclockwise $(i-1)$-rotation and we join its ends to $s_k+j$. 
 We name the resulting path $\mathbf{P_k^{s_{k-1}+i}}$, and we define \[\mathscr{P}_k:=\{\mathbf{P_k^{s_{k-1}+1}},\dots,\mathbf{P_k^{s_k}}\}.\]

 Each path $\mathbf{P_k^{s_{k-1}+i}}$ has the same first and last vertex as its lightface version $P_k^{s_{k-1}+i}$, but spans not only level $k$ (except for one vertex), but every higher level too. 

 It is worth noting that, for each $i\in [p_k]$, the path $\mathbf{P_k^{s_{k-1}+i}}\in\mathscr{P}_k$
 has the following properties:

\begin{itemize}
    \item [\textbf{(P1)}] $\mathbf{P_k^{s_{k-1}+i}}$ has exactly one edge in level $k$ of each type in $im(\phi_k)\cup\{c_k\}$, which appears in no other path of $\mathscr{P}_k$ --- in particular, this implies that each such edge is in exactly one path in~$\mathscr{P}_k$. 
    \item [\textbf{(P2)}] $\mathbf{P_k^{s_{k-1}+i}}$ has exactly two edges in level $k$ of each other type, which it inherits from $P_k^{s_{k-1}+i}$ --- in particular, this implies that each such edge is in exactly two paths of $\mathscr{P}_k$.
    \item [\textbf{(P3)}]  For each $j\in [n-s_k]$,  vertex
    $s_k+j$ is incident to exactly two edges $e,e'$ of $\mathbf{P_k^{s_{k-1}+i}}$, moreover,  
  with $e\in E(\mathbf{P_k^{s_{k-1}+(i+\phi_k(s_k+j) \text{(mod $p_k$)})}})$ and $e'\in E(\mathbf{P_k^{s_{k-1}+(i-\phi_k(s_k+j) \text{(mod $p_k$)})}})$.
\end{itemize}
These three properties follow directly from the construction of $\mathbf{P_k^{s_{k-1}+i}}$.
We also observe that 
\begin{itemize}
    \item [\textbf{(P4)}] there do not exist homotypical edges in any path of $\mathscr{P}_k$ that have counterclockwise distance $c_k$ or $\phi_k(s_k+j)$ for some $j\in [n-s_k]$, 
    \end{itemize}
    because in each path of $\mathscr{P}_k$, one edge of type $\phi_k(s_k+j)$ has been deleted, to be substituted by a path of length $2$ passing through the vertex $s_k+j$, as described above (or, for $c_k$, already from each cycle of $\mathcal{C}_k$ one edge of type $c_k$ has been removed in order to form $\mathcal{P}_k$).  

\begin{example}[cont]
 Figure \ref{fig:3} depicts a choice for the path $\mathbf{P_1^1}$. It is a modification of the generating path $P_1^1$ of Figure \ref{fig:2}. The edge $\{1,4\}$, of type $\phi_1(18)=3$, has been substituted by the path $(1,18,4)$ that contains the vertex $18$. We invite the reader to verify that, similarly, one edge of each of the types $\phi_1(19)=6$, $\phi_1(20)=5$ and $\phi_1(21)=7$ has been substituted by a path of length $2$ passing from $19$, $20$ and $21$, respectively. The component cliques $K_{17}$ and $K_4$ of $K_{21}$ are delineated in red and blue, respectively.
 
     \begin{figure}[ht]
        \centering
        \caption{The path $\mathbf{P_1^1}$ obtained from the generating path $P_1^1$.}
         \label{fig:3}
        $\displaystyle
            \begin {tikzpicture}[-latex ,auto ,node distance =2cm and 3cm ,semithick , state/.style ={draw, circle}]
    \node[state,scale=0.8] (18) at (9,10) {18};
    \node[state,scale=0.8] (19) at (7.5,9) {19};
    \node[state,scale=0.8] (20) at (9,8) {20};
    \node[state,scale=0.8] (21) at (10.5,9) {21};
    \node[state,scale=0.9] (1) at (8,6.7) {1};
    \node[state,scale=0.9] (2) at (6.5,6.5) {2};
    \node[state,scale=0.9] (3) at (5.4,6) {3};
    \node[state,scale=0.9] (4) at (4.5,5.1) {4};
    \node[state,scale=0.9] (5) at (4.2,4) {5};
    \node[state,scale=0.9] (6) at (4.5,2.9) {6};
    \node[state,scale=0.9] (7) at (5.3,1.9){7};
    \node[state,scale=0.9] (8) at (6.5,1.5) {8};
    \node[state,scale=0.9] (9) at (8,1.2) {9};
    \node[state,scale=0.8] (10) at (9.5,1.2) {10};
    \node[state,scale=0.8] (11) at (11,1.5) {11};
    \node[state,scale=0.8] (12) at (12.1,2){12};
    \node[state,scale=0.8] (13) at (13,2.9) {13};
    \node[state,scale=0.8] (14) at (13.2,4.2) {14};
    \node[state,scale=0.8] (15) at (12.8,5.5) {15};
    \node[state,scale=0.8] (16) at (11.4,6.4) {16};
    \node[state,scale=0.8] (17) at (10,6.7) {17};

     \path[-] (1) edge[bend left=15] node {} (18);
     \path[-] (18) edge[bend right=35] node {} (4);
     \path[-]  (4) edge  node {} (13);
     \path[-] (13) edge[bend left=30]  node {} (21);
     \path[-] (21) edge[bend left=15] node {} (6);
     \path[-] (6) edge[bend right=30]  node {} (2);
     \path[-]  (2) edge[bend left=30]  node {} (20);
     \path[-] (20) edge[bend left=30] node {} (7);
     \path[-]  (7) edge[bend right=30]  node {} (5);
     \path[-] (5) edge[bend right=45]  node {} (19);
     \path[-] (19) edge[bend left=30] node {} (16);
     \path[-] (16) edge[bend right=30]  node {} (15);
     \path[-] (15) edge[bend right=30]  node {} (12); 
     \path[-] (12) edge  node {} (3); 
     \path[-] (3) edge  node {} (10);
     \path[-] (10) edge[bend left=30]  node {} (14);
     \path[-] (14) edge[bend right=30]  node {} (9);
     \path[-] (9) edge[bend left=30]  node {} (11);
     \path[-] (11) edge[bend left=30]  node {} (17);

     \draw[red] (8.7,4) ellipse (5.2cm and 3.3cm);
     \draw[blue] (9,9) ellipse (2cm and 1.5cm);
    \end{tikzpicture}
    $
   
\end{figure}
 \end{example}

We continue by constructing sets of paths \[\mathcal{R}_k:=\{R_k^{s_k+1},\dots,R_k^n, R_k^{n+1}\},\] where $R_k^{s_k+j}$, for $j\in [n-s_k]$, contains all the edges of $K_{p_k}$ of type $\phi_k(s_k+j)$ but one, and $R_k^{n+1}$ contains all the edges of $K_{p_k}$ of type $c_k$ but one. Specifically, we pick the  edges excluded by the paths in $\mathcal{R}_k$ to be from distinct paths in $\mathscr{P}_k$; that is, they obey the property \begin{itemize}
    \item[\textbf{(P5)}] 
    If $e, e'\notin E(\bigcup\mathcal{R}_{k})$ are edges of different types, then $e, e'$ belong to different paths from $\mathcal{P}_k$.
    
    \end{itemize}
 This is possible
because, by our construction, 
each path of $\mathscr{P}_k$ has exactly one edge of each type in $im(\phi_k)\cup\{c_k\}$, and that edge is different from those of all other paths in $\mathscr{P}_k$, due to \textbf{(P1)}. Note that the paths in each family $\mathscr{R}_k$ are pairwise edge-disjoint and each spans $K_{p_k}$.



 Then, for each $j\in[p_k]$, we extend the path $R_k^{s_k+j}\neq R_m^{n+1}:=\mathbf{R_m^{n+1}}$ by its \emph{connecting edge}, namely 
 the edge $e(R_k^{s_k+j},\mathbf{P_{k+1}^{s_k+j}})$ if $k\leq m-1$ and $j\leq p_{k+1}$, or 
  the edge $e(R_k^{s_k+j},R_{k+1}^{s_k+j})$ if $k\leq m-1$ and $j\geq p_{k+1}+1$, or the edge $e(R_m^{s_m+j},\mathbf{B_{m+1}^{s_m+j})}$ if $k=m$. We thus obtain a path $\mathbf{R_k^{s_k+j}}$ ending at a vertex of level $k+1$. We define \[\mathscr{R}_k:=\{ \mathbf{R_k^{s_k+1}},\dots,\mathbf{R_k^n},\mathbf{R_k^{n+1}}\}.\] 

Finally, we define \[\mathscr{Q}:=\{\mathbf{Q^1},\dots,\mathbf{Q^n},\mathbf{Q^{n+1}}\}\] where $\mathbf{Q^\ell}$ is the 
concatenation 
 of all the boldface paths of exponent $\ell$. That is, for each $\ell\in [n]$,~$\mathbf{Q^\ell}$ is the concatenation 
  of all paths $\mathbf{R^\ell_r}$ for each $r\in[m]$ such that $s_r+1\leq\ell$, and of $\mathbf{P^\ell_k}$ if there exists $k\in [m]$ such that $s_{k-1}+1\leq\ell\leq s_k$, or $\mathbf{B^\ell_{m+1}}$ otherwise. Moreover, $\mathbf{Q^{n+1}}$ is the concatenation of all paths $\mathbf{R_k^{n+1}}$.
  
  Intuitively, $\mathscr{Q}$ consists of the paths in $\mathscr{P}_1$, the concatenated pairs of paths of equal exponents in $\mathscr{R}_1\times\mathscr{P}_2$, the concatenated triples of paths of equal exponents in $\mathscr{R}_1\times\mathscr{R}_2\times\mathscr{P}_3$, and so on, as well as the path $\mathbf{Q^{n+1}}$.

 \begin{example}[cont]
 For our example on $K_{21}$, the family $\mathscr{Q}$ consists of all the paths $\mathbf{P_1^j}$, where $j\in\{1,\dots,17\}$, and of all the concatenated pairs $(\mathbf{R_1^{17+j},\mathbf{B_2^{17+j}}})$, where $j\in\{1,\dots,4\}$. In Figure \ref{fig:4} one can see the path $\mathbf{Q^{18}}$ that we obtain for a specific choice of the path $\mathbf{R_1^{18}}$. The subpath entirely contained in $K_4$ is $\mathbf{B_2^{18}}$. The rest of the path is our choice of $\mathbf{R_1^{18}}$; it consists of all of the edges of type $\phi_1(18)=3$ except for $\{15,1\}$, as well as a connecting edge between $15$ (the last vertex of $\mathbf{R_1^{18}}$) and $18$ (the first vertex of $\mathbf{B_2^{18}}$).   
 
     \begin{figure}[ht]
        \centering
        \caption{The path $\mathbf{Q^{18}}$ obtained from a particular choice of $\mathbf{R_1^{18}}$ in our running example.}
        \label{fig:4}
            \begin {tikzpicture}[-latex ,auto ,node distance =2cm and 3cm ,semithick , state/.style ={draw, circle}]
    \node[state,scale=0.8] (18) at (9,10) {18};
    \node[state,scale=0.8] (19) at (7.5,9) {19};
    \node[state,scale=0.8] (20) at (9,8) {20};
    \node[state,scale=0.8] (21) at (10.5,9) {21};
    \node[state,scale=0.9] (1) at (8,6.7) {1};
    \node[state,scale=0.9] (2) at (6.5,6.5) {2};
    \node[state,scale=0.9] (3) at (5.4,6) {3};
    \node[state,scale=0.9] (4) at (4.5,5.1) {4};
    \node[state,scale=0.9] (5) at (4.2,4) {5};
    \node[state,scale=0.9] (6) at (4.5,2.9) {6};
    \node[state,scale=0.9] (7) at (5.3,1.9){7};
    \node[state,scale=0.9] (8) at (6.5,1.5) {8};
    \node[state,scale=0.9] (9) at (8,1.2) {9};
    \node[state,scale=0.8] (10) at (9.5,1.2) {10};
    \node[state,scale=0.8] (11) at (11,1.5) {11};
    \node[state,scale=0.8] (12) at (12.1,2){12};
    \node[state,scale=0.8] (13) at (13,2.9) {13};
    \node[state,scale=0.8] (14) at (13.2,4.2) {14};
    \node[state,scale=0.8] (15) at (12.8,5.5) {15};
    \node[state,scale=0.8] (16) at (11.4,6.4) {16};
    \node[state,scale=0.8] (17) at (10,6.7) {17};

     \path[-] (1) edge[bend left=30] node {} (4);
     \path[-] (4) edge[bend left=30] node {} (7);
     \path[-]  (7) edge[bend left=30]  node {} (10);
     \path[-] (10) edge[bend left=30]  node {} (13);
     \path[-] (13) edge[bend left=30] node {} (16);
     \path[-] (16) edge[bend left=30]  node {} (2);
     \path[-]  (2) edge[bend left=30]  node {} (5);
     \path[-] (5) edge[bend left=30] node {} (8);
     \path[-]  (8) edge[bend left=30]  node {} (11);
     \path[-] (11) edge[bend left=30]  node {} (14);
     \path[-] (14) edge[bend left=30] node {} (17);
     \path[-] (17) edge[bend left=30]  node {} (3);
     \path[-] (3) edge[bend left=30]  node {} (6); 
     \path[-] (6) edge[bend left=30]  node {} (9); 
     \path[-] (9) edge[bend left=30]  node {} (12);
     \path[-] (12) edge[bend left=30]  node {} (15);
     \path[-] (15) edge[bend right=45]  node {} (18);
     \path[-] (18) edge[bend left=30] node {} (19);
     \path[-] (19) edge node {} (21);

     \draw[red] (8.7,4) ellipse (5.2cm and 3.3cm);
     \draw[blue] (9,9) ellipse (2cm and 1.5cm);
    \end{tikzpicture}
    \end{figure}
 \end{example} 

The reader can easily check that $\mathscr{Q}$ is indeed a path system, as its elements are concatenations of pairwise internally disjoint paths. It also has size $n+1$. It  only remains to show that  $\mathscr{Q}$ is weakly separating.  

We begin by noting that $E(\mathbf{P_1^1})\cup\dots\cup E(\mathbf{P_k^{s_k}})$ is the set of all edges that touch one of the first $k$ levels. As the paths $\mathbf{R_1^1},\dots,\mathbf{R_{k-1}^{s_{k-1}}}$ are also contained in the first $k$ levels, we have $E(\mathbf{Q^1})\cup\dots\cup E(\mathbf{Q^{s_k}})=E(\mathbf{P_1^1})\cup\dots\cup E(\mathbf{P_k^{s_k}})$, so $E(\mathbf{Q^1})\cup\dots\cup E(\mathbf{Q^{s_k}})$ is also the set of all edges that touch one of the first $k$ levels.  
 
 Therefore, edges that touch lower levels are separated from edges that touch higher levels. 
  Also, edges in level $m+1$ are weakly separated from each other by $\mathbf{Q^{s_m+1}},\cdots,\mathbf{Q^n}$. On the other hand, all the edges that touch a level $k\in [m]$ are contained in the paths of $\mathscr{Q}_k:=\mathscr{P}_k\cup\mathscr{R}_k$, which have distinct exponents. Hence, any path $P\in\mathscr{Q}_k$ extends to a different path $P'\in\mathscr{Q}$. 
 In conclusion, to prove that $\mathscr{Q}$ is a weakly separating path system of $K_n$, it suffices to show that, for each $k\in [m]$, $\mathscr{Q}_k$ weakly separates from each other the edges that touch level $k$.  

Let $e,e'$ be two distinct edges that touch level $k$ but are not contained in it. Then the only paths that may contain them both are those of $\mathscr{P}_k$, as each path of $\mathscr{R}_k$ contains at most one such edge, namely its connecting edge. 
 If there are no paths containing both $e$ and $e'$, the two edges are separated from each other, so suppose that there is a path $\mathbf{P_k^{s_{k-1}+i}}$ containing both $e$ and $e'$. 

If $e$ and $e'$ have a common greater end, say $s_k+j$, then by \textbf{(P3)} the only other path of $\mathscr{P}_k$ that contains, say, $e$ is $\mathbf{P_k^{s_{k-1}+(i-\phi_k(s_k+j) \text{(mod $p_k$)})}}$, and the only other path of $\mathscr{P}_k$ that contains $e'$ is $\mathbf{P_k^{s_{k-1}+(i+\phi_k(s_k+j) \text{(mod $p_k$)})}}$. But these are distinct paths, as $2\phi_k(s_k+j)\neq 0 \text{ (mod $p_k$)}$. So, $e$ and $e'$ are    separated from each other.

If $e$ and $e'$ have distinct greater ends, say $s_k+j_1$ and $s_k+j_2$, then again by \textbf{(P3)} the only other path of $\mathscr{P}_k$ that contains $e$ is either $\mathbf{P_k^{s_{k-1}+(i+\phi_k(s_k+j_1) \text{(mod $p_k$)})}}$ or $\mathbf{P_k^{s_{k-1}+(i-\phi_k(s_k+j_1) \text{(mod $p_k$)})}}$,
and the only other path of $\mathscr{P}_k$ that contains $e'$ is either $\mathbf{P_k^{s_{k-1}+(i+\phi_k(s_k+j_2) \text{(mod $p_k$)})}}$ or $\mathbf{P_k^{s_{k-1}+(i-\phi_k(s_k+j_2) \text{(mod $p_k$)})}}$. Again, as $\phi_k(s_k+j_1)$ and $\phi_k(s_k+j_2)$ are distinct elements of $\{1,\dots,\frac{p_k-1}{2}\}$, these two paths, whichever they are among the two listed cases for each of $e$ and $e'$, are distinct, separating $e$ and $e'$ from each other. 

We deduce that 
\begin{itemize}
    \item[\textbf{(C1)}] the edges that touch level~$k$ but are not contained in it are strongly separated from each other by $\mathscr{Q}_k$. 
\end{itemize}

Now, let $e$ and $e'$ both be contained in level~$k$. If, say, $e$ is in one of the paths in $\mathscr{R}_k$ and $e'$ is not in the same path, then of course $e$ is separated from $e'$ (and $e'$, as long as it is in at least two paths of $\mathscr{Q}_k$, is separated from $e$, because either one of these two paths is in $\mathscr{R}_k$, hence edge-disjoint from the path of $\mathscr{R}_k$ that contains $e$, or they are both in $\mathscr{P}_k$, hence separate $e'$ from all the edges in level~$k$). So, let us assume that, if $e$ is in some path of $\mathscr{R}_k$, then $e'$ is also in the same path. We distinguish three cases, which cover all possible situations, modulo switching the names of $e$ and $e'$.

Case 1.
The edges $e$ and $e'$ are both in the same path of $\mathscr{R}_k$. But then they are homotypical, and their type is in $im(\phi_k)\cup\{c_k\}$. By \textbf{(P1)}, that means that each of $e$ and $e'$ belongs to a different path of~$\mathscr{P}_k$, so they are separated from each other. 

Case 2. 
The edges 
$e$ and $e'$ are not in any of the paths in $\mathscr{R}_k$ because their types are not in $im(\phi_k)\cup\{c_k\}$. But then, by \textbf{(P2)}, each of them appears in exactly two paths of $\mathscr{P}_k$, corresponding to the two paths in which it appears in $\mathcal{P}_k$, therefore they are separated from each other. 

Case 3. The edge $e$ is the unique edge of some type in $im(\phi_k)\cup\{c_k\}$ that is not in any paths in $\mathscr{R}_k$. If the same holds for $e'$, then, by \textbf{(P5)}, $e$ and $e'$ are still separated from each other by $\mathscr{P}_k$. If instead $e'$ is of a type not in $im(\phi_k)\cup\{c_k\}$, then, by \textbf{(P2)}, $e'$ appears in two paths in $\mathscr{P}_k$, and by \textbf{(P1)}, $e$ appears only in one, so again $e'$ is separated from $e$. 

In conclusion, 
\begin{itemize}
    \item[\textbf{(C2)}] the edges that are contained in level $k$ are weakly separated from each other by $\mathscr{Q}_k$. Moreover, the only edges in level $k$ that are not separated from all other edges in level $k$ are those that are in only one path of $\mathscr{Q}_k$.
\end{itemize}
  

Finally, let $e$ be contained in level $k$ and let $e'$ touch level $k$ but not be contained in it. Let $s_k+j$ be the greater end of $e'$. Let $\mathbf{P_k^{s_{k-1}+i}}$ be a path of $\mathscr{P}_k$ containing $e'$, and suppose that it also contains $e$. Then, by \textbf{(P3)} and without loss of generality, $e'$ is also in $\mathbf{P_k^{s_{k-1}+(i+\phi_k(s_k+j) \text{(mod $p_k$)})}}$. By \textbf{(P4)}, there do not exist homotypical edges in any path of $\mathscr{P}_k$ that have counterclockwise distance $\phi_k(s_k+j)$ in $K_{p_k}$, so $e\notin\mathbf{P_k^{s_{k-1}+(i+\phi_k(s_k+j) \text{(mod $p_k$)})}}$. Hence $e'$ is separated from $e$; in fact, if $e$ is in at least two paths in $\mathscr{P}_k$, or in a path in $\mathscr{R}_k$ whose connecting edge is not $e'$, then it too is separated from $e'$. 

We just proved that
\begin{itemize}
    \item[\textbf{(C3)}] 
    \begin{itemize}
        \item[(a)] the edges that touch level $k$ but do not lie within it are separated 
        by $\mathscr{Q}_k$
        from all edges in level~$k$,
        \item[(b)] for each edge $e$ that lies in level $k$ and is in at least two paths of $\mathscr{Q}_k$, one of the following holds
        \begin{itemize}
        \item[(i)]
        $e$ is separated from the edges that touch level $k$ but are not contained in it, or
        \item[(ii)] $e$  is in a path  $P\in\mathscr{P}_k$ and a path  $R\in \mathscr{R}_k$, and $P$ contains   the  connecting edge of $R$.
        \end{itemize}
    \end{itemize}
\end{itemize}

Together, conclusions $\textbf{(C1)}$, $\textbf{(C2)}$ and $\textbf{(C3)}$ imply that $\mathscr{Q}_k$ indeed weakly separates the edges that touch level $k$, which concludes the proof of Theorem \ref{theorem1}. They also list all the possible cases in which strong separation fails among the edges that touch a fixed level, which will be useful for proving Theorem \ref{theorem2}.  

\section{Proof that $ssp(K_n)\leq n+9$}\label{sec4}

In the previous section, we explained that edges that touch lower levels are separated from edges that touch higher levels by $\mathscr{Q}$. We begin this section by noting that those of the latter that are in at least two paths are also separated from most of the former. Indeed, given an edge $e$ touching level $k$ and contained in two paths $Q, Q'\in\mathscr{Q}$, the intersections of these two paths with each level $k'< k$ are distinct (hence edge-disjoint) paths in $\mathcal{R}_{k'}$, so together $Q$ and $Q'$ separate~$e$ from all of the edges that are contained in any of the lower levels. Therefore, the only possibly offending edges are the connecting edges in $Q$ and~$Q'$, some of which may be present in both paths, so that $e$ is not separated from them. In particular, if $\mathscr{Q}$ is such that, for every $k\in [m]$, the paths in $\mathscr{R}_k$ contain pairwise distinct connecting edges, then every edge of~$K_n$ that is in at least two paths of $\mathscr{Q}$ is separated from all other edges. 

Moreover, in Section \ref{sec3}, when examining the separation properties of $\mathscr{Q}_k$, we took care to note when both edges in a pair under examination were separated from each other. We now invite the reader to revisit our conclusions \textbf{(C1)}, \textbf{(C2)} and \textbf{(C3)} and verify that there are only three possible kinds of \emph{bad} edges, that is, edges that may not be separated  by $\mathscr{Q}_k$ (resp. $\mathscr{B}$) from all of the other edges that touch the same level $k\in [m]$ (resp. $m+1$). Namely:

\begin{itemize}
    \item Perhaps some edges in level $m+1$ are bad, as $\mathscr{B}$ is only a weakly separating path system. 
    \item Edges contained in a single path are bad, 
    as they fail to be separated from all of the other edges in the same path. There are exactly $n-s_k+1$ of these edges touching (in fact contained in) each level $k\in [m]$: they are the edges $e_k^{s_k+j}$ of type $\phi(s_k+j)$ that are excluded by the paths $R_k^{s_k+j}$, respectively, as well as the edge $e_k^{n+1}$ of type $c_k$ that is excluded by $R_k^{n+1}$. We shall call these the $e$-\emph{edges} of the corresponding paths in $\mathcal{R}_k$. This is to say, no path in $\mathcal{R}_k$ contains its $e$-edge, and each edge $e_k^{s_k+j}$ (or $e_k^{n+1}$) uniquely determines its path in $\mathcal{R}_k$.  
    \item At most two edges in each path in $\mathcal{R}_k$ that are not separated from its connecting edge are bad, namely those that appear in the two paths of $\mathscr{P}_k$ that contain the connecting edge. Let us call these edges $f_k^{s_k+j}$ and $g_k^{s_k+j}$, under the 
    convention that the path of $\mathscr{P}_k$ that contains $f_k^{s_k+j}$ and the connecting edge of $R_k^{s_k+j}$ has a  smaller exponent than the path of $\mathscr{P}_k$ that contains $g_k^{s_k+j}$ and the connecting edge of $R_k^{s_k+j}$. We also refer to $f_k^{s_k+j}$ and $g_k^{s_k+j}$ as
    the $f$- and $g$-\emph{edges} of the paths in $\mathcal{R}_k$ that contain them. \\ Note that the choice of the path $R_k^{s_k+j}$ (that is, of the edge $e_k^{s_k+j}$) and of the last vertex of this path uniquely determine $f_k^{s_k+j}$ and $g_k^{s_k+j}$. Moreover, the path $\mathbf{R_m^{n+1}}$ does not have a connecting edge, so it also does not have an $f$- and $g$-edge. Finally, we remark that the $f$- and $g$-edge of each path $R_k^{s_k+j}$ share a vertex. Indeed, for any edge touching level $k$ and incident to $s_k+j$, the two paths of $\mathscr{P}_k$ that contain it differ by a $\phi_k(s_k+j)$-rotation (either clockwise or counterclockwise), so the same is true for the connecting edge in $\mathbf{R_k^{s_k+j}}$. In particular, the $f$- and $g$-edge of $\mathbf{R_k^{s_k+j}}$, each of which belongs to a different one of these two paths of $\mathscr{P}_k$ and is the unique edge of type $\phi_k(s_k+j)$ in it, also differ by a $\phi_k(s_k+j)$ rotation. As they are also both of type $\phi_k(s_k+j)$, they share a vertex. This observation is not used in our proof, and is only stated explicitly so that the alert reader will not be surprised by its absence. 
\end{itemize}

We introduce the notations $E_k^{s_k+j}:=\{e_k^{s_k+1},\dots,e_k^{s_k+j}\}$ and $E_k:=E_k^{n+1}$. The sets $F_k^{s_k+j}$, $F_k$, $G_k^{s_k+j}$ and $G_k$ are defined similarly.    

\begin{example}[cont]
    Let us have another look at Figure \ref{fig:3} and Figure \ref{fig:4}. 
    
    In level $2$, each edge of type $1$ lies in a unique path of $\mathscr{B}$ (hence of $\mathscr{Q}$), so there is an edge of type $2$ from which it is not separated. For example, $\{18,19\}$ is not separated from $\{19,21\}$.

    In level $1$, the edge $\{15,1\}$ is the unique edge of type $3$ not in $R_1^{18}$. That is, $e_1^{18}=\{15,1\}$. As a consequence, this edge is not separated from all of the other edges in the path of $\mathscr{P}_1$ that contains it (which is $\mathbf{P_1^4}$, because the edge $\{12,15\}$ is in $\mathbf{P_1^1}$, as seen in Figure \ref{fig:3}, so the edge $\{15,1\}$ is in $\mathbf{P_1^4}$). Similarly, each of the edges $e_1^{19}$ of type $6$, $e_1^{20}$ of type $5$, $e_1^{21}$ of type $7$, and $e_1^{22}$ of type $1$, no matter how they have been chosen (remember that the choice is arbitrary except that they have to originate from pairwise distinct paths in $\mathscr{P}_1$), is not separated from the other edges in its path in $\mathscr{P}_1$.

    Finally, $R_1^{18}$ has the connecting edge $\{15,18\}$. There are two paths in $\mathscr{P}_1$ that contain the edge $\{15,18\}$: $\mathbf{P_1^{12}}$ and $\mathbf{P_1^{15}}$. The unique edge of type $3$ in $\mathbf{P_1^1}$ is $\{12,15\}$, so the unique edge of type $3$ in $\mathbf{P_1^{12}}$ is $f_1^{18}=\{6,9\}$, and in $\mathbf{P_1^{15}}$ it is $g_1^{18}=\{9,12\}$. In the same way we find the $f$- and $g$-edges of each path in $\mathcal{R}_1$. 

    Suppose that we choose $R_1^{19}$ so that $e_1^{19}=\{1,7\}$. This is a valid choice; for $R_1^{18}$ we have $e_1^{18}=\{15,1\}$, which lies in $\mathbf{P_1^4}$, whereas $\{1,7\}$ lies in $\mathbf{P_1^{8}}$, so $\mathbf{(P5)}$ is not violated. In particular, let us choose $1$ to be the last vertex of $R_1^{19}$, so that the connecting edge of $R_k^{19}$ is $\{1,19\}$. We invite the reader to practice by showing that then $f_1^{19}=\{13,2\}$ and $g_1^{19}=\{7,13\}$. \textit{(Hint: look at Figure \ref{fig:3} and mentally rotate the bottom level to find the two paths in $\mathscr{P}_1$ that contain the edge $\{1,19\}$. What is the unique edge of type $\phi_1(19)=6$ in each?)}  
\end{example}

In order to handle both the bad edges and the instances of coinciding connecting edges, it will help to have some extra wiggle space. It follows easily from \cite{nagura} (see Lemma \ref{lem2} in Appendix \ref{appendix1}) that for $n\geq 3$ there exists a decomposition $n=\sum_{k=1}^mp_k+b$ with each $p_k$ an odd prime so that for the partial sums $s_k$ we have either $n-s_k\leq \frac{p_k}{4}$ and $b\in\{3,\dots,19, 29, 31\}$, or 
$(p_m,b)\in\big\{(19,5), (19,6), (23,7), (31,8)\big\}$. 

 Our plan is to choose the paths of $\mathcal{R}_k$ more carefully, so that in the end all the bad edges in level $k$ will be neatly sorted in at most six linear forests, and no coinciding connecting edges will exist in distinct paths of $\mathscr{R}_k$. For the remainder of this section, for the sake of ease and consistency, let us follow the convention that the last vertex of a path in $\mathcal{R}_k$, out of its two extremal vertices $u$ and $v$, is chosen to be $v$ if and only if $u-v \textbf{ (mod $p_k$)}\in \lceil{\frac{p_k-1}{2}\rceil}$ --- in particular, the choice of an $e$-edge uniquely determines the corresponding oriented path in $\mathcal{R}_k$ and its $f$- and $g$-edges.  

We say that level $k$ is \emph{arranged} if each of $E_k$, $F_k$ and $G_k$ can be decomposed into at most two linear forests. When tackling the weak version of the problem, the only restriction that we imposed on the choice of $\mathcal{R}_k$ was that the edges in $E_k$ be drawn from different paths of $\mathscr{P}_k$, which is property \textbf{(P5)}. Now we demand in addition that 
 \begin{itemize} 
     \item[\textbf{(P6)}] level $k$ is arranged, and
      \item[\textbf{(P7)}] no two of the paths of $\mathcal{R}_k$ have the same connecting edge.
 \end{itemize}


Naturally, level $m+1$ is arranged, as all of the 
sets $E_{m+1}$, $F_{m+1}$ and $G_{m+1}$ are empty. We inductively arrange all the levels of $K_n$, down to the first, as follows. We sequentially choose each oriented path $R_k^{s_k+j}$ and orient each maximal path in $E_k^{s_k+j}$, $F_k^{s_k+j}$ and $G_k^{s_k+j}$ (in short, we \emph{establish} $R_k^{s_k+j}$), according to the following process. The path $R_k^{s_k+1}$ is established arbitrarily. Once the paths $\{R_k^{s_k+1}, \dots, R_k^{s_k+j}\}$ have been established, we choose the path $R_k^{s_k+j+1}$ so that:

\begin{enumerate} [label=(\roman*)]
     \item the path in $\mathscr{P}_k$ that contains $e_k^{s_k+j+1}$ does not contain any edge in $E_k^j$ (this is just \textbf{(P5)});
    \item the $e$-edge (resp. $f$-edge, $g$-edge) of $R_k^{s_k+j+1}$ is not incident to any end of an edge in $E_k^j$ (resp. $F_k^j$, $G_k^j$), unless said end is the last vertex of a maximal path of $E_k^j$ (resp. $F_k^j$, $G_k^j$) (we will see below that this ensures \textbf{(P6)});
    \item if the connecting edge of $R_k^{s_k+j+1}$ has the same greater end as the connecting edge of a previously established path $R_k^{s_k+j'+1}$, $j'<j$, then $R_k^{s_k+j+1}$ may not have the same last vertex as $R_k^{s_k+j'+1}$ (this is just \textbf{(P7)}). 
\end{enumerate}

We then arbitrarily orient each newly formed maximal path in $E_k^{s_k+j+1}$, $F_k^{s_k+j+1}$ and $G_k^{s_k+j+1}$ and repeat. This process continues for as long as we can obey the above three conditions, and then halts. 

Observe that item (ii) ensures that, when the above process cannot continue anymore, the bad edges of each kind that will have been already selected 
form a linear forest. Indeed, each of $E_k^{s_k+1}$, $F_k^{s_k+1}$ and $G_k^{s_k+1}$ has at most one edge, hence is a linear forest. If $E_k^{s_k+j}$, $F_k^{s_k+j}$ and $G_k^{s_k+j}$ are linear forests, then adding to each of them a single edge that is not incident to an internal vertex of any path and does not connect the ends of any single path also produces a linear forest.

We will iterate the entire process up to two times to arrange level $k$. At the end of two iterations, we must have selected a total of $n-s_k+1$ paths, which is the required size of $\mathcal{R}_k$. Let us begin.

Jointly, items (i) and (ii) forbid at most $7j$ options for $R_k^{s_k+j+1}$ (that is, for $e_k^{s_k+j}$). Indeed, by (i), each edge of $E_k^{s_k+j}$ forbids us from choosing $e_k^{s_k+j+1}$ from a unique path of $\mathscr{P}_k$, and each such path contains a unique edge of type $\phi_k(s_k+j+1)$, so $j$ options are forbidden because of that. Moreover, if $F$ is a linear forest consisting of n oriented paths and $F$ has $j$ edges, then $F$ has exactly $j$ vertices that are not last in any of these paths, as can be easily proven by induction. By (ii), each non-last vertex in the linear forest $E_k^{s_k+j}$ forbids as candidates for $e_k^{s_k+j+1}$ the two edges of type $\phi_k(s_k+j+1)$ that are incident to it. Similarly, each of the $j$ non-last vertices in $F_k^{s_k+j}$ and $G_k^{s_k+j}$ forbids two options for $f_k^{s_k+j+1}$ and $g_k^{s_k+j+1}$, respectively, hence for $e_k^{s_k+j+1}$. This sums to $6j$ forbidden options at most, which, together with the $j$ options from item (i), give us at most $7j$ forbidden options.        

On the other hand, item (iii) forbids up to four other options. Indeed, each of the vertices of level $k+1$ may be the first vertex of up to five paths in $\mathscr{Q}_{k+1}$: one in $\mathscr{P}_{k+1}$, as the paths in $\mathscr{P}_{k+1}$ have pairwise distinct first vertices, and at most four in $\mathscr{R}_{k+1}$, because each end of a path in $\mathscr{R}_k$ is also an end of its $e$-edge, and, as level $k+1$ is arranged by induction, the $e$-edges of level $k+1$ form at most two linear forests, meaning that at most four of them are incident to any given vertex. This implies that for each vertex of level $k+1$ there may be up to five paths in $\mathcal{R}_k$ that connect to it with their connecting edges. So, when a path of $\mathcal{R}_k$ is being established, up to four vertices of level $k$ are forbidden from being its last vertex for this reason.  

An exception to the above is if $k=m$, in which case item (iii) is irrelevant; the paths of $\mathscr{B}$ have pairwise distinct first vertices, so the connecting edges that are incident to their first vertices are necessarily distinct. 

Therefore, letting $A$ denote the index at which the first iteration stops, we are able to establish in total $A\geq\min\{\lceil\frac{p_k-4}{7}\rceil+1, n-s_k+1\}$ (or $A\geq\min\{\lceil\frac{p_k}{7}\rceil, n-s_k+1\}$ if $k=m$) paths in $\mathcal{R}_k$ so that the properties \textbf{(P5)} and $\textbf{(P7)}$ are satisfied and each of $E_k^{s_k+A}$, $F_k^{s_k+A}$ and $G_k^{s_k+A}$ is a linear forest. 

If $A\leq n-s_k$, we select $R_k^{s_k+A+1}$ arbitrarily (but obeying $\mathbf{(P5)}$) and we begin the second iteration of the process. This time, in each step we apply item (ii) only to the sets $E_k^j\setminus E_k^{s_k+A}$, $F_k^j\setminus F_k^{s_k+A}$ and $G_k^j\setminus G_k^{s_k+A}$. On the other hand, items (i) and (iii) still concern the entire set of previously established paths, including those from the first iteration. Reasoning exactly as for the first iteration, we see that, in the establishment of the path $R_k^{s_k+j+1}$ with $j\geq A$, there are at most $j$ options that are forbidden due to item (i), and at most $6(j-A)$ options that are forbidden due to item (ii). Moreover, there are up to four options that are forbidden by item (iii), unless  $k=m$, in which case there are none. Hence we can add $B\geq\min\{\lceil\frac{p_k-A-4}{7}\rceil, n-s_k-A+1\}$ (or $B\geq\min\{\lceil\frac{p_k-A}{7}\rceil, n-s_k-A+1\}$ if $k=m$) more paths in $\mathcal{R}_k$ so that the properties \textbf{(P5)} and $\textbf{(P7)}$ are satisfied and each of $E_k^{s_k+A+B}\setminus E_k^{s_k+A}$, $F_k^{s_k+A+B}\setminus F_k^{s_k+A}$ and $G_k^{s_k+A+B}\setminus G_k^{s_k+A}$ is a linear forest.

For $p_k\geq 67$, one can show by routine calculation and the fact $n-s_k\leq\frac{p_k}{4}$ that $A+B\geq n-s_k+1$, so this second iteration is indeed the final one, establishing the required $n-s_k+1$ paths in $\mathcal{R}_k$ in such a way that $\textbf{(P5)}$, $\textbf{(P6)}$ and $\textbf{(P7)}$ are obeyed. Almost all of the smaller cases, including the special cases with $(p_m,b)=(19,5)$ or $(31,8)$, can be checked individually to the same effect, that is,~that $A+B\geq b+1$ (because for $p_k\leq 66$, if we have $n-s_k\leq \frac{p_k}{4}$, then $n-s_k\leq 17$, so $k=m$ and $n-s_k=b$). 

For the cases with $(p_m,b)=(19,6)$ or $(23,7)$, the process yields only six and seven paths, respectively. However, in both cases  slight tweaks yield the desired additional path. For $(19,6)$, after establishing $\lceil\frac{19}{7}\rceil=3$ paths during the first iteration, there are at most five options for $e_m^{s_m+4}$ (hence for $R_m^{s_m+4}$) that would cause $E_m^{s_m+4}$ to cease being a linear forest. Indeed, if in this case we exceptionally allow for the new edge to be incident to the first vertex of a maximal path (as long as it does not create a cycle), then the greatest possible number of forbidden options for $e_m^{s_m+4}$ is five (rather than six), only potentially achieved when the edges of $E_m^{s_m+3}$ form a path (Figure \ref{fig5}). Similarly, there are at most five options for each of $f_m^{s_m+4}$ and $g_m^{s_m+4}$ that are forbidden, each giving a forbidden option for $e_m^{s_m+4}$. As $e_m^{s_m+4}$ must also originate from a different path of $\mathscr{P}_m$ than the three edges of $E_m^{s_m+3}$, we have a total of at most eighteen forbidden options for $e_m^{s_m+4}$, hence for $R_m^{s_m+4}$. As there are nineteen options in total, this allows us to establish one more path in the first iteration of the process. The second iteration is performed normally, obtaining $\lceil\frac{19-4}{7}\rceil=3$ more paths, for a total of seven. For $(23,7)$, we run the first iteration normally and establish $\lceil\frac{23}{7}\rceil=4$ paths, then an additional four are established in the second iteration of the process by using the exact same method as for the first iteration of $(19,6)$. That is, we establish $\lceil\frac{23-4}{7}\rceil=3$ paths normally, then we notice that, out of the twenty-three options for $R_m^{s_m+8}$, only eighteen ar forbidden.

\begin{figure}[h!]
\centering
\caption{Given the linear forest $E_k^{s_k+3}$ (or $F_k^{s_k+3}$, or $G_k^{s_k+3}$) with three edges, there exist at most five edges of type $\phi(s_k+4)$ that can be added to it to produce a graph that is not a linear forest.}
\label{fig5}
\begin {tikzpicture}[-latex ,auto ,node distance =2cm and 3cm ,semithick , state/.style ={draw, circle}]
    \node[state,scale=0.8] (1) at (9,3) {};
    \node[state,scale=0.8] (2) at (11,1) {};
    \node[state,scale=0.8] (3) at (13,1) {};
    \node[state,scale=0.8] (4) at (15,3) {};
    \node[scale=0.8] (5) at (11.5,2.5) {};
    \node[scale=0.8] (6) at (12.5,2.5) {};
    \node[scale=0.8] (7) at (11.5,-0.5) {};
    \node[scale=0.8] (8) at (12.5,-0.5) {};

    \path[-] (1) edge[red] node {} (4);
    \path[-] (1) edge[] node {} (2);
    \path[-] (2) edge[] node {} (3);
    \path[-] (3) edge[] node {} (4);
    \path[-] (2) edge[red] node {} (5);
    \path[-] (3) edge[red] node {} (6);
    \path[-] (2) edge[red] node {} (7);
    \path[-] (3) edge[red] node {} (8);
    \end{tikzpicture}
\end{figure}


Finally, for each $k\in [m]$, we extend each path in $\mathcal{R}_k$ by its connecting edge to create the family $\mathscr{R}_k$. Subsequently, we use the paths of $\mathscr{P}_k$ and $\mathscr{R}_k$ for each $k\in [m]$, and of $\mathscr{B}$, to obtain again a weakly separating path system $\mathscr{Q}$ of $K_n$ exactly as in Section \ref{sec3}, albeit one that has the additional properties $\textbf{(P6)}$ and $\textbf{(P7)}$. In particular, every edge of $K_n$ that is in at least two paths of $\mathscr{Q}$ is separated from all other edges.

\begin{example} [cont]
    Let us run the above process on our example of $K_{21}$. First of all, note that we can keep the same partition into $K_{17}$ and $K_4$ that we had in Section \ref{sec3}, because $21-17\leq\frac{17}{4}$. Let us also choose $R_1^{18}$ to be the same and with the same last vertex $15$, as this choice is arbitrary. We need to complete a total of five steps in order to establish the paths $R_1^{18},\dots,R_1^{22}$. 

    Before we begin, let us describe exactly what it is that we will do. This will be useful to the reader who wishes to replicate our process. As we explained in our analysis, at each step $j$ for $j\in\{1,\dots,4\}$ we must firstly choose a candidate edge $e_1^{17+j}$ of type $\phi_1(17+j)$. This will determine the path $R_1^{17+j}$ and a last vertex for this path. This then determines the edges $f_1^{17+j}$ and $g_1^{17+j}$. At this point, we check whether one of the items (i), (ii) or (iii) is violated. If so, we choose a different edge $e_1^{17+j}$ and repeat. Otherwise, we orient the newly formed maximal paths of $E_1^{17+j}$, $F_1^{17+j}$ and $G_1^{17+j}$ and we proceed to the next step. The preceding analysis proves that, in just two iterations of this process, we will manage to complete a sufficient number of steps to define all five paths of $\mathcal{R}_1$. However, it does not offer a deterministic algorithm by which to make valid choices at each step, except for the obvious: exhaustive search. At step $j+1$, we arbitrarily select an edge of type $\phi(s_k+j+1)$ as a candidate for $e_k^{s_k+j+1}$, which in turn determines also the edges $f_k^{s_k+j+1}$ and $g_k^{s_k+j+1}$. If our choice is valid (that is, it does not lead to violation of items (i), (ii) or (iii)), we arbitrarily reorient the newly formed maximal paths in $E_k^{s_k+j+1}$, $F_k^{s_k+j+1}$ and $G_k^{s_k+j+1}$ and proceed to the next step. Otherwise, we arbitrarily select a new candidate edge that we have not yet tried, and repeat until we succeed.  

    In what follows, we will write each maximal path of $E_1^j$, $G_1^j$ or $F_1^j$, for each $j\in \{1,\dots,17\}$, in the form $[u_1,\dots,u_r]$, where each two consecutive vertices are adjacent in the path and $u_r$ is the last vertex of the path. 

    To begin, we have $E_1^{18}=e_1^{18}=[1,15]$, $F_1^{18}f_1^{18}=[6,9]$, and $G_1^{18}=g_1^{18}=[12,9]$. 

    Could we continue with the $R_1^{19}$ of our exercise? At first glance, everything looks good. As we have already mentioned, \textbf{(P5)} is not violated. Moreover, it is impossible to create a graph that is not a linear forest with only two edges. Finally, we do not need to worry about whether the last vertices of $R_1^{18}$ and $R_1^{19}$ coincide (although they do not, as the former is $15$ and the latter is $1$), because we are at level $k=1$, and we have $m=1$ in our example, so $k=m$. However, the last vertex of the path $e_1^{18}=[1,15]$ is $15$, and $e_1^{19}=\{1,7\}$ is incident to $e_1^{18}$ at a vertex that is not its last! Does this matter? Could we not simply backtrack and say that the last vertex of the path $e_1^{18}$ is $1$? Yes, but the point of our analysis is to eschew the need for fine control and ad hoc fixes. Once more we state that the last vertices of maximal paths are not carefully selected, but are arbitrary. So, we reject this choice of $R_1^{19}$. Instead, we keep looking, and we come upon a choice with $e_1^{19}=[8,2]$, last vertex $2$, $f_1^{19}=[3,14]$, and $g_1^{19}=[8,14]$, which can be readily checked to be valid.

    Fortunately, after only one iteration of the process rather than two, we arrived at the results listed in Table \ref{tableprocess}. Of course, the reader who attempts this experiment is expected to obtain a different (but correct) outcome depending on the choices made (or sampled). It is possible that two iterations will be required.
\end{example}

     \begin{table}[ht]
        \centering
        \caption{A valid selection of the paths $R_1^{18},\dots,R_1^{22}$ in our example.}
        \label{tableprocess}
        \hspace*{-0.6cm}\begin{tabular}{c|c|c|c|c|c|c|c}
           $e_1^{17+j}$ & $\mathscr{P}_1$-path with $e_1^{17+j}$ & last vertex of $R_1^{17+j}$ & $f_1^{17+j}$ & $g_1^{17+j}$ & $E_1^{17+j}$ & $F_1^{17+j}$ & $G_1^{17+j}$\\
            \hline
            $\{15,1\}$ & $P_1^4$ & $15$ & $\{6,9\}$ & $\{9,12\}$ & $[1,15]$ & $[6,9]$ & $[12,9]$\\
            \hline
            $\{2,8\}$ & $P_1^9$ & $2$ & $\{14,3\}$ & $\{8,14\}$ & $[1,15]$ & $[6,9]$ & $[12,9]$\\
            & & & & & $[8,2]$ & $[3,14]$ & $[8,14]$\\
            \hline
            & & & & & $[1,15]$ & $[6,9]$ & $[12,9]$\\
            $\{4,9\}$ & $P_1^{13}$ & $4$ & $\{11,16\}$ & $\{16,4\}$ & $[8,2]$ & $[3,14]$ & $[8,14]$\\
            & & & & & $[4,9]$ & $[11,16]$ & $[4,16]$\\
            \hline 
             & & & & & $[1,15]$ & $[6,9]$ & $[12,9,2]$\\
             $\{5,12\}$ & $P_1^3$ & $5$ & $\{9,16\}$ & $\{2,9\}$ & $[8,2]$ & $[3,14]$ & $[8,14]$\\
            & & & & & $[4,9]$ & $[11,16,9]$ & $[4,16]$\\ 
            & & & & & $[5,12]$ & &\\
            \hline
             & & & & & $[1,15]$ & $[6,9]$ & $[12,9,2]$\\
            $\{9,10\}$ & $P_1^{12}$ & $9$ & - & - & $[8,2]$ & $[3,14]$ & $[8,14]$\\
            & & & & & $[4,9,10]$ & $[11,16,9]$ & $[4,16]$\\ 
            & & & & & $[5,12]$ & &
        \end{tabular}
        \label{tab:placeholder}
    \end{table}

It remains to separate the bad edges from all of the other edges. We point out two helpful facts. Firstly, the $e$-edges in the entire graph are strongly separated from each other by $\mathscr{Q}$, and need only to be 
separated from all of the other edges. Indeed, for $e$-edges that touch distinct levels, the one in the lower level is separated from the one in the higher level as is the case for all edges, whereas the path in $\mathscr{Q}$ that contains the one in the higher level avoids the one in the lower level, as the latter is not in any $\mathscr{R}$-path. Moreover, $e$-edges from the same level are strongly separated from each other by $\mathscr{Q}$. 
Secondly, the $f$- and $g$-edges in each level only need to be 
separated from the connecting edges that touch that level.

Let us begin by constructing some additional paths that separate the $e$-, $f$- and $g$-edges from all the other edges of $K_n$. Note that we will not modify any of the paths of $\mathscr{Q}$, but we will come up with more paths, so that all of the other cases of edge separation that are already achieved by~$\mathscr{Q}$ will be preserved. 

For each linear forest made up of the $f$-edges of level $k\leq m-1$ obtained during one of the two iterations of our process, we can connect its separate paths by edges between their ends to form a single path. We may do so arbitrarily, as each $f$-edge in level $k$ is already separated from every other edge in its level by $\mathscr{Q}_k$. In this manner, we form in each level $k\leq m$ two paths $F_k^{\alpha}$ and $F_k^{\beta}$, each corresponding to a different iteration of our process, whose union contains all of the $f$-edges of level $k$. Even though only one iteration of the process may have been required to arrange level $k$, we can always partition the $f$-edges (of which, recall, there are $n-s_k+1\geq 2$) arbitrarily into two non-empty subsets, so that
$F_k^{\alpha}$ and $F_k^{\beta}$ are both non-empty.

Note that for each $k\leq m-1$ there are at least $p_k-2(n-s_k+1)+2=p_k-2(n-s_k)$ 
edges between the first vertex of $F_{k+1}^{\alpha}$ and the vertices of level $k$ that are not internal in $F_k^{\alpha}$. The $f$-edges of level $k$ are not separated from at most $n-s_k+1$ of them, namely the connecting edges in paths of $\mathscr{R}_k$. As $p_k-3(n-s_k)-1\geq\frac{p_k}{4}-1>0$, we can use one of the rest, if needed together with an edge connecting its end in level $k$ with the last vertex of $F_k^{\alpha}$, to join $F_k^{\alpha}$ and $F_{k+1}^{\alpha}$ into a single path. In the exact same manner we can connect the last vertex of $F_k^{\beta}$ to the first vertex of $F_{k+1}^{\beta}$ by a path of length at most $2$. Doing this for every $k\leq m-1$ from top to bottom creates two paths $\mathbf{F^{A}}$ and $\mathbf{F^{B}}$ that live in the first $m$ levels and separate all the $f$-edges from all of the connecting edges. The exact same tactic yields two paths $\mathbf{G^{A}}$ and $\mathbf{G^{B}}$ that live in the first $m$ levels and accomplish the same for the $g$-edges.

We cannot achieve the same level of frugality for the $e$-edges. However, it is easy to check that the edges of a linear forest $F$ inside a complete graph can be separated from the rest of the edges by at most two paths whose intersection is precisely $F$. That is, there exist paths $E^{\alpha}$, $E^{\beta}$, $E^{\gamma}$ and $E^{\delta}$ that separate the $e$-edges from all of the other edges. Let us say, without loss of generality, that $E^{\alpha}$ and $E^{\beta}$ contain the edges obtained in the first iteration of our process, whereas $E^{\gamma}$ and $E^{\delta}$ contain those that arose in the second iteration. 

Now, let us turn our attention to the edges of level $m+1$, the only potentially bad edges remaining. If $b=12$ or $15$, then $\mathscr{B}$ is a strongly separating path system (see Lemma \ref{small_paths}, Appendix \ref{appendix2}) and we are done; $\mathscr{Q}$ already separates the edges of level $m+1$ from each other, and each edge in level $m+1$ is in two paths of $\mathscr{Q}$, so it is separated from all the edges in lower levels. Otherwise, let $c_{m+1}$ be the common type of the edges contained in level $m+1$ that are not separated from all other edges by $\mathscr{B}$; the type $c_{m+1}$ exists by $\textbf{(W3)}$ and is given for each possible $b$ in Lemma \ref{small_paths} Appendix \ref{appendix2}. The edges of type $c_{m+1}$ in level $m+1$ form $gcd(b,c_{m+1})\leq 2$ vertex-disjoint cycles. We can decompose these at most two cycles into two vertex-disjoint paths $P_E$, $P_E'$ and two vertex-disjoint edges $e_R$, $e_R'$. We will join these four paths to other paths that we have already described. The reader should keep in mind that, whenever two paths in a path system are joined into a single path (perhaps with the addition of a connecting edge), the edges in this new path may not be separated from each other anymore. So, for each of the joinings that we do below, we explain why the edges in the resulting path remain separated from each other.   

For each of the paths $E^{\alpha}$, $E^{\beta}$, $E^{\gamma}$ and $E^{\delta}$ (which are all disjoint from $P_E$ and $P_E'$, as they do not intersect level $m+1$), we connect one of its ends to the first vertex of $P_E$ and the other to the first vertex of $P_E'$ (for $E^{\alpha}$ and $E^{\gamma}$) or we connect one of its ends to the last vertex of $P_E$ and the other to the last vertex of $P_E'$ (for $E^{\beta}$ and $E^{\delta}$), to obtain four paths $\mathbf{E^A}$, $\mathbf{E^B}$, $\mathbf{E^\Gamma}$ and $\mathbf{E^\Delta}$. The intersection of these four new paths is exactly $P_E\cup P_E'$, so they separate the edges of $P_E$ and $P_E'$ from all other edges of~$K_n$. Moreover, the $e$-edges of $K_n$ remain separated from the edges of $P_E$ and $P_E'$ by $\mathscr{Q}$, and as the edges connecting $E^{\alpha}$ and $E^{\beta}$ to $P_E$ and $P_E'$ are distinct (and similarly for $E^{\gamma}$ and $E^{\delta}$), the $e$-edges also remain separated from those. Finally, each of the four edges used to join the four $E$-paths with $P_E$ is in two paths of $\mathscr{Q}$, so it is separated from all other edges of $K_n$.  

We also connect with edges the first and the last vertex of $\mathbf{R^{n+1}}$ to $e_R$ and $e_R'$, respectively, to form a path $\mathbf{R}$. The joining edge that connects the first vertex of $\mathbf{R^{n+1}}$ to $e_R$ must not also be the connecting edge in the path of $\mathscr{Q}$ that contains $e_R$, otherwise $e_R$ will not be separated from it. As there are two choices for this joining edge, one for each end of $e_R$, this is possible. Similarly, the joining edge between the last vertex of $\mathbf{R^{n+1}}$ and $e_R'$ is not the connecting edge in the path of $\mathscr{Q}$ that contains $e_R'$. Now, the path $\mathbf{R}$ separates the edges $e_R$ and $e_R'$ from all the other edges in level $m+1$. Each of $e_R$ and $e_R'$ is also in a path of $\mathscr{Q}$, and $\mathbf{R^{n+1}}$ is edge-disjoint from the $\mathscr{R}$-tail of every other path in $\mathscr{Q}$ (because, for each $k\in [m]$, every two paths of $\mathscr{R}_k$ are edge-disjoint), so $e_R$ and $e_R'$ are now also separated from every edge that touches one of the first $m$ levels, so in total they are separated from all other edges. The edges of $\mathbf{R^{n+1}}$ remain separated from $e_R$ and $e_R'$, and from the two joining edges, by $\mathscr{Q}\setminus\mathbf{R^{n+1}}$. Finally, each of the two joining edges is in two paths of $\mathscr{Q}$, so it is separated from all other edges of $K_n$.

 Let $\mathscr{Q}':=(\mathscr{Q}\setminus \mathbf{R^{n+1}})\cup \mathbf{R}$. We conclude that the set \[
 \mathscr{Q}'\cup\{\mathbf{E^A}, \mathbf{E^B}, \mathbf{E^\Gamma}, \mathbf{E^\Delta}, \mathbf{F^A}, \mathbf{F^B}, \mathbf{G^A}, \mathbf{G^B}\}\] strongly separates $K_n$. The proof is complete.

\bibliographystyle{plain}
\bibliography{main}

\begin{appendices} 

\section{Sum decompositions}\label{appendix1}

\begin{lemma}\label{lem1}
    Every natural number $n\geq 3$ can be written in the form $n=\sum_{k=1}^m p_k+b$, where each $p_k$ is an odd prime number such that $n-\sum_{i=1}^kp_i\leq \frac{p_k-3}{2}$ and $b\in\{3,\dots,19\}$.
\end{lemma}

\begin{proof}
    By the main result of \cite{nagura}, for every $n\geq 30$ there exists a prime number $p$ such that $\frac{5}{6}n\leq p\leq n$. Therefore, for $n\geq 33$, there exists a prime number $p$ such that $\frac{5}{6}(n-3)\leq p\leq n-3$. This implies that
    \begin{equation}\label{eq1}
        n-p\geq 3.
    \end{equation}
    Moreover, for $n\geq 21$, we have $\frac{5}{6}(n-3)\geq \frac{2}{3}n+1$, so for $n\geq 33$ we have $\frac{2}{3}n+1\leq p$, which implies 
    \begin{equation}\label{eq2}
        n-p\leq \frac{p-3}{2}.
    \end{equation}
 We also verify that for $n\in\{20,\dots,32\}$ there is a prime number $p$ such that $n$ and $p$ satisfy Equations \ref{eq1} and \ref{eq2} (see Table \ref{tab:placeholder}). The lemma follows by finding such primes $p_1$ for $n$, $p_2$ for $n-p_1$, etc. until we are left with $b\in\{3,\dots,19\}$.  
\end{proof}

\begin{lemma}\label{lem2}
    Every natural number $n\geq 3$ can be written in the form $n=\sum_{k=1}^m p_k+b$, where each $p_k$ is an odd prime number such that one of the following holds:
    \begin{enumerate}[(a)]
        \item $n-\sum_{i=1}^kp_i\leq \frac{p_k}{4}$ and $b\in\{3,\dots,19,29,31\}$, or
        \item  $(p_m,b)\in\big\{(19,5), (19,6), (23,7),(31,8)\big\}$.
    \end{enumerate}
    
\end{lemma}

\begin{proof}
    The proof is almost identical to that of Lemma \ref{lem1}. However, $n-p\leq\frac{p}{4}\Leftrightarrow p\geq\frac{4}{5}n$, and the inequality $\frac{5}{6}(n-3)\geq\frac{4}{5}n$ is true only for $n\geq 75$. This means that we have to check our statements $n-p\ge 3$ and $n-p\le \frac p4$ separately for all natural numbers $20\leq n\leq 74$ (see Table \ref{tab:placeholder}). In doing so for all values of $n$, we find six exceptions: two ($29$ and $31$) that are themselves allowed values for $b$ as listed in item (a) of the lemma, and four ($24$, $25$, $30$ and $39$) 
    for which we can choose appropriate values of $b$ as listed in item~(b) of the lemma. 
\end{proof}

\begin{table}[ht]
    \centering
    \caption{For each natural number $n$ in $\{20,\dots,74\}$ we give a prime number $p$ such that $n-p\geq 3$, and $n-p\leq\frac{p}{4}\leq\frac{p-3}{2}$ (black entries) or $n-p\leq\frac{p-3}{2}$ (red entries).}
    \label{tab:placeholder}
    \begin{tabular}{c|c}
       n & p \\
       \hline
        20 & 17\\
       \hline
        21 & 17\\
       \hline
        22 & 19\\
       \hline
        23 & 19\\
       \hline
        {\color{red}24} & {\color{red}19}\\
       \hline
        {\color{red}25} & {\color{red}19}\\
       \hline
        26 & 23\\
       \hline
        27 & 23\\
       \hline
        28 & 23\\
       \hline
        {\color{red}29} & {\color{red}23}\\
       \hline
        {\color{red} 30} & {\color{red} 23}\\
       \hline
        {\color{red}31} & {\color{red}23}\\
       \hline
        32 & 29\\
       \hline
        33 & 29\\
    \end{tabular}
    \hspace*{1cm}\begin{tabular}{c|c}
         n & p \\
         \hline
        34 & 31\\
       \hline
        35 & 31\\
       \hline
        36 & 31\\
       \hline
        37 & 31\\
        \hline
        38 & 31\\
         \hline
        {\color{red}39} & {\color{red}31}\\
       \hline
        40 & 37\\
        \hline
       41 & 37\\
       \hline
        42 & 37\\
       \hline
        43 & 37\\
       \hline
        44 & 41\\
       \hline
        45 & 41\\
       \hline
        46 & 43\\
       \hline
        47 & 43\\
    \end{tabular}
    \hspace*{1cm}\begin{tabular}{c|c}
        n & p \\
        \hline
        48 & 43\\
       \hline
        49 & 43\\
       \hline
        50 & 47\\
       \hline
        51 & 47\\
       \hline
        52 & 47\\
       \hline
        53 & 47\\
       \hline
        54 & 47\\
       \hline
        55 & 47\\
       \hline
        56 & 53\\
        \hline
        57 & 53\\
       \hline
        58 & 53\\
       \hline
        59 & 53\\
        \hline
        60 & 53\\
       \hline
        61 & 53\\
    \end{tabular}
    \hspace*{1cm}\begin{tabular}{c|c}
        n & p \\
        \hline
        62 & 59\\
        \hline
        63 & 59\\
       \hline
        64 & 61\\
       \hline
        65 & 61\\
       \hline
        66 & 61\\
        \hline
        67 & 61\\
       \hline
        68 & 61\\
       \hline
        69 & 61\\
        \hline
       70 & 67\\
       \hline
        71 & 67\\
       \hline
        72 & 67\\
       \hline
        73 & 67\\
        \hline
       74 & 71\\
    \end{tabular}
\end{table}

\section{Some facts about weakly separating path or cycle systems}\label{appendix2}

In this appendix we include for completeness two components of \cite{wickes} that are important in our proof. 

Firstly, we claim that every complete graph on a prime number of vertices greater than or equal to $5$ contains a generating cycle. This is our property \textbf{(W1)}. Although this is not explicitly mentioned by Wickes in~\cite{wickes}, it can be easily surmised from her proof of her Lemma 3.4.

\begin{lemma}[Lemma 3.4, \cite{wickes}]
    There exists a generating path for $K_n$ whenever $n$ is an odd prime.
\end{lemma}

Let us briefly review that proof. Given such a clique $K_p$, Wickes labels its vertices with the elements of the field $\mathbb{F}_p$ and takes a primitive root $g$ of $\mathbb{F}_p$, that is, a generator of the multiplicative subgroup $U(\mathbb{F}_p)$. She then considers the walk $P:=[1,1+g,\dots,1+\sum_{i=1}^{p-2}g^i]$, and proceeds to show that this is a generating path in $K_p$. Specifically, it has two edges of each type except for type $1$, of which it has one edge. Moreover, the two edges of type $g^j$ (or $p-g^j$, whichever is in $[\frac{p-1}{2}]$) have as their counterclockwise distance some number contained in the set $$\Big\{\pm\frac{g^{\frac{p-1}{2}+j+1}-g^j}{g-1}\text{ (mod p)},\pm\frac{g^{\frac{p-1}{2}+j}-g^{j+1}}{g-1}\text{ (mod p)}\Big\}$$ and the counterclockwise distances of distinct homotypical pairs are distinct.

With very minor modifications, we can obtain, instead of a generating path, a generating cycle.

\begin{lemma}\label{generating_cycle}
    There exists a generating cycle for $K_n$ whenever $n\geq 5$ is a prime number.
\end{lemma}

\begin{proof}

To obtain a generating cycle, we simply add to $P$ the edge $$e=\Big\{1+\sum_{i=1}^{p-2}g^i,1+\sum_{i=1}^{p-1}g^i\Big\}.$$ The resulting walk $C$ is indeed a cycle, as $1+\sum_{i=1}^{p-1}g^i=\frac{g^p-1}{g-1}=1$. Moreover, it has two edges of each type. Indeed, Wickes \cite{wickes} shows that $P$ has two edges of each type except for type $1$, of which it has one. In $C$, the edge $e$ is also of type $g^{p-1}=1$. Finally, that the counterclockwise distance of the two edges of type $1$ is distinct from the counterclockwise distances of all other homotypical pairs can be proved exactly as done for other homotypical pairs in Lemma 3.4 of \cite{wickes} without any modifications.

We remark that this proof does not work for $p=3$ because then the ``cycle" that we obtain has length $2$, but our graph is simple.

\end{proof}

Secondly, we claim that, for each $n\in\{3,\dots,19\}$, there exists in $K_n$ a path whose set of images under rotation is a weakly separating path system of $K_n$, and the edges that violate strong separation, if they exist, are all of a single type, for which we use the term ``special type" below. Examples of such paths for each $n\in\{3,\dots,19\}$ appeared first in \cite{wickes}, and we list them again here. For each of the paths, we note if the weakly separating path system that results from its rotation is in fact strongly separating, or else which is the unique type that pertains to edges that are not separated from all other edges.\vspace{5mm}

\begin{lemma}\label{small_paths}
    For each $n\in\{3,\dots,19\}$, there exists in $K_n$ a path whose images under rotation form a weakly separating path system of $K_n$, and the edges that violate strong separation, if any, are of a single type.
\end{lemma}

\begin{proof}

$P(3)=[1,2]$, with special type $1$\\
$P(4)=[1,2,4]$, with special type $1$\\
$P(5)=[1,3,2,5]$, with special type $1$\\
$P(6)=[1,5,4,3,6]$, with special type $2$\\
$P(7)=[1,2,3,5,7,4]$, with special type $3$\\
$P(8)=[1,3,5,2,6,7,8]$, with special type $3$\\
$P(9)=[1,5,9,3,4,6,8,2]$, with special type $1$\\
$P(10)=[1,4,7,6,5,9,3,8,10]$, with special type $2$\\
$P(11)=[1,3,5,10,4,11,7,8,9,6]$, with special type $3$\\
$P(12)=[1,2,11,9,10,3,7,4,8,6,12,5]$, strongly separating\\
$P(13)=[1,3,4,13,11,6,10,7,12,5,8,9]$, with special type $6$\\
$P(14)=[1,3,6,9,10,11,2,7,13,5,12,8,4]$, with special type $2$\\
$P(15)=[1,14,15,5,10,3,12,6,9,13,2,4,11,8,7]$, strongly separating\\
$P(16)=[1,11,13,15,14,3,8,12,16,9,2,10,7,4,5]$, with special type $6$\\
$P(17)=[1,3,5,16,10,11,12,9,6,15,7,14,4,17,13,8]$, with special type $5$\\
$P(18)=[1,15,10,5,13,3,12,9,6,7,8,2,14,16,18,11,4]$, with special type $4$\\
$P(19)=[1,3,5,18,12,11,10,13,16,7,17,6,14,9,4,19,15,8]$, with special type $7$\vspace{5mm}

Verifying the special types (and their uniqueness) is not too time-consuming. In most of the above paths, the edges in each homotypical pair are consecutive, so their counterclockwise distance is their type. Moreover, in the above cases in which $n$ is even, each edge of type $\frac{n}{2}$ can be rotated by $\frac{n}{2}$ to coincide with itself. That is, all of the counterclockwise distances of homotypical pairs of edges are pairwise distinct (where each edge of type $\frac{n}{2}$ is considered to be in a homotypical pair with itself rotated by $\frac{n}{2}$). As proved by Wickes (Lemma 3.3, \cite{wickes}), this implies that the edges in these homotypical pairs are separated from all other edges. Hence, the only edges that fail to be separated from all other edges are those of a type that appears only once in the corresponding listed path (and is not $\frac{n}{2}$). For each of the paths above there is at most one type for which this occurs, namely the one listed. 

The only paths listed above in which the edges of homotypical pairs are not conveniently placed next to each other are $P(5)$, $P(9)$, $P(12)$, $P(13)$, $P(15)$ and $P(16)$, which require careful inspection.

\end{proof}

\end{appendices}

\end{document}